\RequirePackage{fix-cm}
\documentclass[smallcondensed]{svjour3}       % onecolumn (second format)
\smartqed  % flush right qed marks, e.g. at end of proof
\usepackage{graphicx}
\usepackage[utf8]{inputenc}
\usepackage[english]{babel}
\usepackage{amsmath}
\usepackage{amssymb,amsfonts}
\usepackage{mathrsfs}
\usepackage{url}
\newcommand{\N}{\mathbb{N}}
\newcommand{\Z}{\mathbb{Z}}

\newcommand{\R}{\mathbb{R}}

 \journalname{Bollettino dell'Unione Matematica Italiana}

\begin{document}
\title{ On the interplay between hypergeometric series, Fourier-Legendre expansions and Euler sums
%\thanks{Grants or other notes
%about the article that should go on the front page should be
%placed here. General acknowledgments should be placed at the end of the article.}
}
\subtitle{}

\titlerunning{Hypergeometric series, FL-expansions and Euler sums}        % if too long for running head

\author{Marco Cantarini         \and
        Jacopo D'Aurizio %etc.
}

%\authorrunning{Short form of author list} % if too long for running head

\institute{Marco Cantarini \at
              Universit\`a di Perugia\\
Orcid 0000-0003-2310-0939\\
              Tel.: +39-3296369936\\
 \email{cantarini\_m@libero.it}  
         \\    
        %\\
%             \emph{Present address:} of F. Author  %  if needed
        \and
         Jacopo D'Aurizio \at
 Universit\`a di Pisa\\
 Orcid 0000-0002-4965-748X\\
              Tel.: +39-3479835357\\
 \email{jacopo.daurizio@gmail.com}  
}

\date{Received: date / Accepted: date}
% The correct dates will be entered by the editor

\maketitle

\begin{abstract}
In this work we continue the investigation, started in \cite{Interplay}, about the interplay between hypergeometric functions and Fourier-Legendre ($\textrm{FL}$) series expansions. 
In the section ``\emph{Hypergeometric series related to $\pi,\pi^2$ and the lemniscate constant}'', through the FL-expansion of $\left[x(1-x)\right]^\mu$ (with $\mu+1\in\frac{1}{4}\N$) 
we prove that all the hypergeometric series 
$$
 \sum_{n\geq 0}\frac{(-1)^n(4n+1)}{p(n)}\left[\frac{1}{4^n}\binom{2n}{n}\right]^3,\quad \sum_{n\geq 0}\frac{(4n+1)}{p(n)}\left[\frac{1}{4^n}\binom{2n}{n}\right]^4,$$ 
 $$\quad \sum_{n\geq 0}\frac{(4n+1)}{p(n)^2}\left[\frac{1}{4^n}\binom{2n}{n}\right]^4,\; \sum_{n\geq 0}\frac{1}{p(n)}\left[\frac{1}{4^n}\binom{2n}{n}\right]^3,\; \sum_{n\geq 0}\frac{1}{p(n)}\left[\frac{1}{4^n}\binom{2n}{n}\right]^2 $$
return rational multiples of $\frac{1}{\pi},\frac{1}{\pi^2}$ or the lemniscate constant, as soon as $p(x)$ is a polynomial fulfilling suitable symmetry constraints. Additionally, by computing 
the FL-expansions of $\frac{\log x}{\sqrt{x}}$ and related functions, we show that in many cases the hypergeometric $\phantom{}_{p+1} F_{p}(\ldots , z)$ function evaluated
at $z=\pm 1$ can be converted into a combination of Euler sums. In particular we perform an explicit evaluation of
$$ \sum_{n\geq 0}\frac{1}{(2n+1)^2}\left[\frac{1}{4^n}\binom{2n}{n}\right]^2,\quad \sum_{n\geq 0}\frac{1}{(2n+1)^3}\left[\frac{1}{4^n}\binom{2n}{n}\right]^2. $$
In the section ``\emph{Twisted hypergeometric series}'' we show that the conversion of some $\phantom{}_{p+1} F_{p}(\ldots,\pm 1)$ values into combinations of Euler sums, 
driven by FL-expansions, applies equally well to some \emph{twisted} hypergeometric series, i.e. series of the form $\sum_{n\geq 0} a_n b_n$ where 
$a_n$ is a Stirling number of the first kind and $\sum_{n\geq 0}b_n z^n = \phantom{}_{p+1} F_{p}(\ldots;z)$.
\keywords{Hypergeometric functions \and Fourier - Legendre expansions \and Euler sums}
% \PACS{PACS code1 \and PACS code2 \and more}
 \subclass{33C20 \and 30B50}
\end{abstract}

\section{List of abbreviations}
\begin{itemize}
 %\item $\text{Beta}$ : properties of Euler Beta function;
 \item $K$ : Catalan's constant, $\sum_{n\geq 0}\frac{(-1)^n}{(2n+1)^2}$;
 \item $\text{IBP}$ : integration by parts;
 %\item $\text{Li}_2$ : functional identities for the dilogarithm function;
 %\item $\text{Li}_3$ : functional identities for the trilogarithm function;
 \item $\text{FL}$ : Fourier-Legendre expansion;
 %\item $\text{FC}_1$ : Expansion by shifted Chebyshev polynomials of the first kind;
 %\item $\text{FC}_2$ : Expansion by shifted Chebyshev polynomials of the second kind;
 \item $\text{PFD}$ : partial fraction decomposition;
 \item $\text{Taylor}$ : Taylor series at the origin;
 \item $(a)_k$: rising Pochhammer symbol, $\Gamma(a+k)/\Gamma(a)$.
\end{itemize}
In order to simplify the notation we will denote the complete elliptic integrals of the first and second
kind as $K(x)$ and $E(x)$, regarding $x$ as the elliptic modulus, such that $K(x)=\text{{\bf K}}(\sqrt{x})$ and 
$E(x)=\text{{\bf E}}(\sqrt{x})$. Additionally, we will employ the simplified notation for hypergeometric functions:
$$ {}_{3}F_{2}\left[ \begin{matrix} \frac{2a+1}{4}, \frac{2a+3}{4}, \frac{2b+1}{2} \\ c, 
 \frac{2 d + 1}{2} \end{matrix} \ \Bigg| \ 1 \right] = \phantom{}_3 F_2\left(\tfrac{2a+1}{4},\tfrac{2a+3}{4},\tfrac{2b+1}{2}; c,\tfrac{2d+1}{2}; 1\right).$$
Unless differently stated, we will also assume that the variables involved in our series range over the interval $(0,1)$.

\section[Hypergeometric series related to pi and the lemniscate constant]{Hypergeometric series related to $\pi,\pi^2$ and the lemniscate constant}
\label{intro}
In our previous works (\cite{Parbelos},\cite{Interplay}) we stressed many consequences of the following principle:

 \textbf{Remark}. Many hypergeometric functions have a Fourier-Legendre expansion whose coefficients still have a hypergeometric structure.
 In particular, given a hypergeometric series $\phantom{}_p F_q(\ldots;\pm 1)=\sum_{n\geq 0}a_n = \int_{0}^{1} f(x)\,g(x)\,dx$ and assuming that both $f$ an $g$ have 
 manageable FL-expansions, the original series can be converted into an equivalent one. In many practical cases, the transformed series 
 can be directly computed from Euler sums with low ($\leq 5$) weight.
$$\,$$
Such approach turned out to be pretty effective in computing hypergeometric series involving $\left[\frac{1}{4^n}\binom{2n}{n}\right]^2$ 
of $\left[\frac{1}{64^n}\binom{4n}{2n}\binom{2n}{n}\right]$. Starting with the identity
$$ K(x)\stackrel{L^2(0,1)}{=}\sum_{n\geq 0}\frac{2}{(2n+1)} P_n(2x-1) $$
we proved that many series involving central binomial coefficients or their squares can be computed in explicit terms, most of the times 
depending on $\pi, K, \log(2),$ $\log(1+\sqrt{2})$ and $\Gamma\left(\frac{1}{4}\right)$. By considering the moments 
$$ \int_{0}^{1} x^\eta K(x)\,dx, \qquad \int_{0}^{1} x^\eta E(x)\,dx $$
and similar integrals associated to \emph{generalized} elliptic integrals, we proved that the FL-machinery is able to produce 
representations for $\pi$ in terms of ratios of $\phantom{}_3 F_2$ functions. Such representations have been usually proved 
through the contiguity relations for hypergeometric functions \cite{Vidunas}. The current section is about the application of the same idea 
to hypergeometric functions of the $\frac{\left[x(1-x)\right]^n}{\sqrt{x(1-x)}}$ kind, with $n\in\mathbb{N}$. In particular,
the starting point of our current investigation is the identity
\begin{equation}\label{uno} \frac{1}{\sqrt{x(1-x)}}\stackrel{L^2(0,1)}{=}\pi \sum_{n\geq 0}(4n+1)\left[\frac{1}{4^n}\binom{2n}{n}\right]^2 P_{2n}(2x-1) \end{equation}
which is a simple consequence of Rodrigues' formula for Legendre polynomials. By Bonnet's recursion formula, if the 
FL-expansion of $f(x)$ is given by $\sum_{n\geq 0} c_n P_n(2x-1)$, then the FL-expansion of $g(x)=x\cdot f(x)$ is given by 
$\sum_{n\geq 0}d_n P_n(2x-1)$, where for any $n\geq 1$
$$ d_n = \frac{n}{2(2n-1)} c_{n-1} + \frac{1}{2}c_n + \frac{n+1}{2(2n+3)} c_{n+1}. $$
Additionally, the symmetry of Legendre polynomials ensures that the FL-expansion of $f(1-x)$ is simply given by $\sum_{n\geq 0} c_n(-1)^n P_n(2x-1)$. 
In particular it is straightforward to exploit (\ref{uno}) to compute the FL-expansion of $\sqrt{\frac{x}{1-x}}$ and $\sqrt{\frac{1-x}{x}}$, then the FL-expansion 
of $\sqrt{x(1-x)}$:
\begin{equation}\label{due} \sqrt{x(1-x)}\stackrel{L^2(0,1)}{=} \frac{\pi}{8}\sum_{n\geq 0}\frac{4n+1}{(n+1)(1-2n)}\left[\frac{1}{4^n}\binom{2n}{n}\right]^2 P_{2n}(2x-1). \end{equation}
Since $\frac{1}{4^n}\binom{2n}{n}\sim\frac{1}{\sqrt{\pi n}}$ as $n\to +\infty$, the series $\sum_{n\geq 0}\frac{4n+1}{(n+1)(1-2n)}\left[\frac{1}{4^n}\binom{2n}{n}\right]^2$
is absolutely convergent and the equality above holds as a \emph{pointwise} equality for any $x\in[0,1]$, too. Since $P_{2n}(0)=(-1)^n\left[\frac{1}{4^n}\binom{2n}{n}\right]$, 
the evaluation of (\ref{due}) at $x=\frac{1}{2}$ leads to:
\begin{equation}\label{tre} \frac{1}{\pi} = \sum_{n\geq 0}\frac{(4n+1)(-1)^n}{4(n+1)(1-2n)}\left[\frac{1}{4^n}\binom{2n}{n}\right]^3\end{equation}
and the application of Parseval's identity leads to:
\begin{equation}\label{quattro} \frac{1}{\pi^2} = \sum_{n\geq 0}\frac{3(4n+1)}{32(n+1)^2(1-2n)^2}\left[\frac{1}{4^n}\binom{2n}{n}\right]^4,\end{equation}
while by computing the inner product between $\sqrt{\frac{x}{1-x}}$ and $\sqrt{\frac{1-x}{x}}$ we have:
\begin{equation}\label{quattrobis} \frac{16}{\pi^2} = \sum_{n\geq 0}\frac{8n^2+8n+1}{(n+1)^2}\left[\frac{1}{4^n}\binom{2n}{n}\right]^4,\end{equation}
which are Ramanujan-like formulas for $\frac{1}{\pi}$ and $\frac{1}{\pi^2}$. 
The RHS of (\ref{quattrobis}) depends on $\phantom{}_4 F_3\left(\tfrac{1}{2},\tfrac{1}{2},\tfrac{1}{2},\tfrac{1}{2};\begin{smallmatrix}1, & 1, & 1 \\ 1, & 1, & 2 \\ 1, & 2, & 2\end{smallmatrix}; 1\right)$.
The formula above can be seen as a consequence of the contiguity relations for $\phantom{}_4 F_3$ functions; nevertheless, it does not seem to be recognized 
by computer algebra systems. We may notice that for (moderately) large values of $n$ the main term in the
RHS of (\ref{tre}) behaves like $ \frac{C}{n^{5/2}}$ and the main term of the RHS of (\ref{quattro}) behaves like $\frac{C}{n^5}$. In terms of the 
contiguity relations for hypergeometric functions, (\ref{tre}) is equivalent to
$$ \frac{16}{\pi} = 4\cdot\phantom{}_3 F_2\left(-\tfrac{1}{2},\tfrac{1}{2},\tfrac{1}{2};1,2;-1\right)+\phantom{}_3 F_2\left(\tfrac{1}{2},\tfrac{3}{2},\tfrac{3}{2};2,3;-1\right)$$
and (\ref{quattro}) is equivalent to
$$ \frac{512}{3\pi^2} = 16\cdot\phantom{}_4 F_3\left(-\tfrac{1}{2},-\tfrac{1}{2},\tfrac{1}{2},\tfrac{1}{2};1,2,2;1\right)+\phantom{}_4 F_3\left(\tfrac{1}{2},\tfrac{1}{2},\tfrac{3}{2},\tfrac{3}{2};2,3,3;1\right).$$
On the other hand, the procedure which allowed us to derive the FL-expansion of $\sqrt{x(1-x)}$ from the FL-expansion of $\frac{1}{\sqrt{x(1-x)}}$ can be applied again.
It leads to the FL-expansions of $\left[x(1-x)\right]^{3/2}$ and $\left[x(1-x)\right]^{5/2}$:

$$ \left[x(1-x)\right]^{3/2} = \frac{9\pi}{256}\sum_{n\geq 0}\frac{(4n+1)}{(n+1)_2 \left(\frac{1}{2}-n\right)_2}\left[\frac{1}{4^n}\binom{2n}{n}\right]^2 P_{2n}(2x-1), $$
$$ \left[x(1-x)\right]^{5/2} = \frac{225\pi}{4096}\sum_{n\geq 0}\frac{(4n+1)}{(n+1)_3 \left(\frac{1}{2}-n\right)_3}\left[\frac{1}{4^n}\binom{2n}{n}\right]^2 P_{2n}(2x-1). $$
Except for the extra term $\left[\frac{1}{4^n}\binom{2n}{n}\right]^2$, these expansions are essentially identical to the FL-expansions for the generalized 
elliptic integrals we studied in our previous work. The evaluation at $x=\frac{1}{2}$ produces:
\begin{equation}\label{cinque} \frac{1}{\pi} = \sum_{n\geq 0}\frac{9(4n+1)(-1)^n}{8(n+1)(n+2)(1-2n)(3-2n)}\left[\frac{1}{4^n}\binom{2n}{n}\right]^3 \end{equation}
\begin{equation}\label{sei} \frac{1}{\pi} = \sum_{n\geq 0}\frac{225(4n+1)(-1)^n}{16(n+1)(n+2)(n+3)(1-2n)(3-2n)(5-2n)}\left[\frac{1}{4^n}\binom{2n}{n}\right]^3 \end{equation}
where, for (moderately) large values of $n$, the main term of the RHS of (\ref{cinque}) behaves like $\frac{1}{n^{9/2}}$ and the main term of the RHS of (\ref{sei}) behaves like $\frac{1}{n^{13/2}}$.
(\ref{cinque}) and (\ref{sei}) are equivalent to the contiguity relations
$$ \frac{32}{3\pi} = 2\cdot\phantom{}_3 F_2\left(-\tfrac{3}{2},\tfrac{1}{2},\tfrac{1}{2};1,3;-1\right)+\phantom{}_3 F_2\left(-\tfrac{1}{2},\tfrac{3}{2},\tfrac{3}{2};2,4;-1\right)$$
$$ \frac{256}{5\pi} = 8\cdot\phantom{}_3 F_2\left(-\tfrac{5}{2},\tfrac{1}{2},\tfrac{1}{2};1,4;-1\right)+5\cdot\phantom{}_3 F_2\left(-\tfrac{3}{2},\tfrac{3}{2},\tfrac{3}{2};2,5;-1\right)$$
and Parseval's identity ensures
\begin{equation}\label{sette} \frac{1}{\pi^2}=\frac{2835}{16384}\sum_{n\geq 0}\frac{(4n+1)}{(n+1)_2^2\left(\frac{1}{2}-n\right)_2^2}\left[\frac{1}{4^n}\binom{2n}{n}\right]^4, \end{equation}
\begin{equation}\label{otto} \frac{1}{\pi^2}=\frac{35083125}{4194304}\sum_{n\geq 0}\frac{(4n+1)}{(n+1)_3^2 \left(\frac{1}{2}-n\right)_3^2}\left[\frac{1}{4^n}\binom{2n}{n}\right]^4, \end{equation}
where, for (moderately) large values of $n$, the main term of the RHS of (\ref{sette}) behaves like $\frac{1}{n^{9}}$ and the main term of the RHS of (\ref{otto}) behaves like $\frac{1}{n^{13}}$.
(\ref{sette}) and (\ref{otto}) are equivalent to the contiguity relations
$$ \frac{16384}{315\pi^2} = 4\cdot\phantom{}_4 F_3\left(-\tfrac{3}{2},-\tfrac{3}{2},\tfrac{1}{2},\tfrac{1}{2};1,3,3;1\right)+\phantom{}_4 F_3\left(-\tfrac{1}{2},-\tfrac{1}{2},\tfrac{3}{2},\tfrac{3}{2};2,4,4;1\right),$$
$$ \frac{16777216}{17325\pi^2} = 64\cdot\phantom{}_4 F_3\left(-\tfrac{5}{2},-\tfrac{5}{2},\tfrac{1}{2},\tfrac{1}{2};1,4,4;1\right)+25\cdot\phantom{}_4 F_3\left(-\tfrac{3}{2},-\tfrac{3}{2},\tfrac{3}{2},\tfrac{3}{2};2,5,5;1\right).$$
\begin{theorem}\label{FLuno} By induction on $k$ we get that for any odd $k\in\mathbb{N}$ the following identity holds pointwise over $[0,1]$:

$$\left[x(1-x)\right]^{\frac{k}{2}} = \frac{\pi(k!!)^2}{8^k 2^{\frac{k+1}{2}}}\sum_{n\geq 0}\frac{(4n+1)\left[\frac{1}{4^n}\binom{2n}{n}\right]^2}{(n+1)_{\frac{k+1}{2}}\cdot \left(\frac{1}{2}-n\right)_{\frac{k+1}{2}}}\,P_{2n}(2x-1).$$

\end{theorem}
\begin{corollary} If $p(x)\in\mathbb{Z}[x]$ is a polynomial with even degree, with simple roots belonging to the set $\left\{\ldots,-3,-2,-1,\frac{1}{2},\frac{3}{2},\frac{5}{2},\ldots\right\}$, 
and such that $p(x)=p\left(-\frac{1}{2}-x\right)$,
$$ \sum_{n\geq 0}(-1)^n\frac{4n+1}{p(n)}\left[\frac{1}{4^n}\binom{2n}{n}\right]^3 $$
is a rational multiple of $\frac{1}{\pi}$ and 
$$ \sum_{n\geq 0}\frac{4n+1}{p(n)^2}\left[\frac{1}{4^n}\binom{2n}{n}\right]^4,\qquad \sum_{n\geq 0}\frac{4n+1}{p(n)}\left[\frac{1}{4^n}\binom{2n}{n}\right]^4 $$
are rational multiples of $\frac{1}{\pi^2}$.
\end{corollary}
\begin{proof} If $p(n)$ fulfills the mentioned constraints, $\frac{1}{p(n)}$ can be written as a linear combination with rational coefficients of 
$ \frac{4n+1}{(n+1)(1-2n)}, \frac{4n+1}{(n+1)(n+2)(1-2n)(3-2n)}$ etcetera. The claim is so a straightforward consequence of Theorem $\ref{FLuno}$.\qed 
\end{proof}
By considering the inner product between different FL-expansions provided by Theorem $\ref{FLuno}$ we also get that
\begin{corollary} If $p(x)\in\mathbb{Z}[x]$ is a monic polynomial with even degree, with simple \emph{or double} roots belonging 
to the set $\left\{\ldots,-3,-2,-1,\frac{1}{2},\frac{3}{2},\frac{5}{2},\ldots\right\}$, 
and such that $p(x)=p\left(-\frac{1}{2}-x\right)$, then 
$$ \sum_{n\geq 0}\frac{4n+1}{p(n)}\left[\frac{1}{4^n}\binom{2n}{n}\right]^4 $$
is a rational multiple of $\frac{1}{\pi^2}$.
\end{corollary}
The evaluation of the identity of Theorem \ref{FLuno} at $x=\frac{1}{2}$ leads to
\begin{equation}\label{nove}\sum_{n\geq 0}\frac{(-1)^n(4n+1)\left[\frac{1}{4^n}\binom{2n}{n}\right]^3}{(n+1)_k\left(\frac{1}{2}-n\right)_k} =\frac{2}{\pi}\cdot \frac{16^k k!^2}{(2k)!^2}\end{equation}
and if we apply Parseval's identity to the statement of Theorem \ref{FLuno} we get:
\begin{equation}\label{dieci} \sum_{n\geq 0}\frac{(4n+1)\left[\frac{1}{4^n}\binom{2n}{n}\right]^4}{(n+1)_k^2 \left(\frac{1}{2}-n\right)_k^2} = \frac{2^{12k}(k!)^4}{\pi^2 k(4k)!(2k)!^2}.  \end{equation}
Another function of the $\left[x(1-x)\right]^\eta$ kind has a simple FL-expansion which is related to central binomial coefficients:
\begin{equation}\label{minusfourth}\left[x(1-x)\right]^{-1/4} = \frac{4\pi\sqrt{\pi}}{\Gamma\left(\frac{1}{4}\right)^2}\sum_{n\geq 0}\left[\frac{1}{4^n}\binom{2n}{n}\right]P_{2n}(2x-1) \end{equation}
the evaluation at $x=\frac{1}{2}$ and Parseval's identity lead to:
$$ \frac{2}{\pi}K(-1)=\sum_{n\geq 0}(-1)^n\left[\frac{1}{4^n}\binom{2n}{n}\right]^2 = \frac{\Gamma\left(\frac{1}{4}\right)^2}{(2\pi)^{3/2}},$$ $$\qquad \frac{2}{\pi}\int_{0}^{1}K(x^4)\,dx=\sum_{n\geq 0}\frac{1}{4n+1}\left[\frac{1}{4^n}\binom{2n}{n}\right]^2 = \frac{\Gamma\left(\frac{1}{4}\right)^4}{16\pi^2} $$
while by combining (\ref{uno}) and (\ref{minusfourth}) we get:
$$ \sum_{n\geq 0}\left[\frac{1}{4^n}\binom{2n}{n}\right]^3 = \frac{\Gamma\left(\frac{1}{4}\right)^4}{4\pi^3}.$$
This provides a furher proof of $K\left(\frac{1}{2}\right)=\frac{1}{4\sqrt{\pi}}\Gamma\left(\frac{1}{4}\right)^2$, not directly relying on Clausen's formula, 
giving the equality between $\phantom{}_2 F_1\left(a,b;a+b+\frac{1}{2};z\right)^2$ and\\ $\phantom{}_3 F_2\left(2a,2b,a+b;a+b+\frac{1}{2},2a+2b;z\right)$.
Given the generating function for\\ Legendre polynomials, the identity $\frac{1}{\sqrt{1-x}}\stackrel{L^2(0,1)}{=}2\sum_{n\geq 0}P_n(2x-1)$ is straightforward.
If we combine the FL-expansions of $\frac{1}{\sqrt{1-x}}$ and $\frac{1}{\sqrt[4]{x(1-x)}}$ with Theorem (\ref{FLuno}) we obtain
the following extension of Bauer's and Levrie's results (\cite{Bauer},\cite{Levrie}):
\begin{theorem}\label{FLdue} For any $k\in\mathbb{N}^+$ we have
 $$\sum_{n\geq 0}\frac{\left[\frac{1}{4^n}\binom{2n}{n}\right]^3}{(n+1)_k\left(\frac{1}{2}-n\right)_k} = \frac{\Gamma\left(\frac{1}{4}\right)^2}{2\pi\sqrt{2}}\cdot \frac{\Gamma\left(k+\tfrac{1}{4}\right)}{\Gamma\left(k+\tfrac{1}{2}\right)^2\,\Gamma\left(k+\tfrac{3}{4}\right)},$$
 $$\sum_{n\geq 0}\frac{\left[\frac{1}{4^n}\binom{2n}{n}\right]^2}{(n+1)_k\left(\frac{1}{2}-n\right)_k} = \frac{2^{8k-3}\Gamma(k)^2}{\pi\,\Gamma(4k)}.$$
Additionally, if $p(x)\in\mathbb{Z}[x]$ is a monic polynomial with even degree, with simple roots belonging to the set $\left\{\ldots,-3,-2,-1,\frac{1}{2},\frac{3}{2},\frac{5}{2},\ldots\right\}$, 
and such that $p(x)=p\left(-\frac{1}{2}-x\right)$, then 
$$ \sum_{n\geq 0}\frac{1}{p(n)}\left[\frac{1}{4^n}\binom{2n}{n}\right]^3 $$
is a rational multiple of $\frac{\Gamma\left(\frac{1}{4}\right)^4}{8\pi^3}$ and 
$$ \sum_{n\geq 0}\frac{1}{p(n)}\left[\frac{1}{4^n}\binom{2n}{n}\right]^2 $$
is a rational multiple of $\frac{1}{\pi}$. 
\end{theorem}
We may also notice that by combining the FL-expansions of $\frac{1}{\sqrt{x(1-x)}}$\\ and $\frac{1}{\sqrt{2-x}}$ we have
\begin{eqnarray*} \int_{0}^{1}\frac{dx}{\sqrt{x(1-x)(2-x)}}&\stackrel{x\mapsto x^2}{=}&\int_{0}^{1}\frac{2\,dx}{\sqrt{(1-x^2)(2-x^2)}}=\sqrt{2}\,K\left(\tfrac{1}{2}\right)\\ &\stackrel{\text{FL}}{=}& 2\pi \sum_{n\geq 0}(\sqrt{2}-1)^{4n+1}\left[\frac{1}{4^n}\binom{2n}{n}\right]^2\\
&\stackrel{\text{Taylor}}{=}& 4(\sqrt{2}-1)\,K\left((\sqrt{2}-1)^4\right),\end{eqnarray*}
which can be deduced from $K(x)=\frac{\pi}{2\,\text{AGM}\left(1,\sqrt{1-x}\right)}$. Such identity relates $\mathbf{K}(k_1),$ $\mathbf{K}(k_4)$ and the lemniscate constant, 
where $k_j=\lambda^*(j)$ are the elliptic singular values, algebraic numbers fulfilling $\frac{\mathbf{K}'(k_r)}{\mathbf{K}(k_r)}=\sqrt{r}$.\\
If we combine Theorem (\ref{FLuno}) with the FL-expansion for the logarithm
$$ \log(x)\stackrel{L^2(0,1)}{=}-1+\sum_{n\geq 1}(-1)^{n+1}\left(\frac{1}{n}+\frac{1}{n+1}\right)P_n(2x-1) $$
we get that for any odd $k\in\mathbb{N}$:
\begin{eqnarray*} &&\frac{\Gamma\left(1+\frac{k}{2}\right)^2}{(k+1)!}\left[H_{k/2}-H_{k+1}\right] = \int_{0}^{1}\log(x)\left[x(1-x)\right]^{k/2}\,dx\\  & = & -\frac{\pi(k!!)}{8^k\left(\frac{k+1}{2}\right)!}-\frac{\pi(k!!)^2}{8^k\cdot 2^{\frac{k+1}{2}}}\sum_{n\geq 1}\frac{\left(\frac{1}{2n}+\frac{1}{2n+1}\right)\left[\frac{1}{4^n}\binom{2n}{n}\right]^2}{(n+1)_{\frac{k+1}{2}}\left(\frac{1}{2}-n\right)_{\frac{k+1}{2}}}
\end{eqnarray*}
and by letting $k=(2j-1)$:
$$ \frac{(2j)!}{4^j}\sum_{n\geq 1}\frac{\left(\frac{1}{2n}+\frac{1}{2n+1}\right)\left[\frac{1}{4^n}\binom{2n}{n}\right]^2}{(n+1)_{j}\left(\frac{1}{2}-n\right)_{j}} = 2\log(2)-1-\sum_{m=1}^{2j}\frac{(-1)^{m+1}}{m}$$
which is an acceleration formula relating $\phantom{}_4 F_3\left(1,1,\tfrac{3}{2},\tfrac{3}{2}-j;2,2,2+j;1\right)$ and\\ $\phantom{}_3 F_2\left(\tfrac{3}{2},\tfrac{3}{2},\tfrac{3}{2}-j;\tfrac{5}{2},j+2;1\right)$.
By Rodrigues' formula the FL-expansion of $\sqrt[4]{x(1-x)}$ is given by
$$ \sqrt[4]{x(1-x)} = \frac{\Gamma\left(\frac{1}{4}\right)^2}{4\sqrt{\pi}}\sum_{n\geq 0}\frac{4n+1}{(4n+3)(1-4n)}\left[\frac{1}{4^n}\binom{2n}{n}\right]P_{2n}(2x-1) $$
and by considering the inner product between the FL-expansions of\\ $\left[x(1-x)\right]^{\pm 1/4}$\\ we immediately get:
$$ \frac{1}{\pi}=\sum_{n\geq 0}\frac{1}{(4n+3)(1-4n)}\left[\frac{1}{4^n}\binom{2n}{n}\right]^2, $$
which is equivalent to $\phantom{}_3 F_2\left(-\tfrac{1}{4},\tfrac{1}{2},\tfrac{1}{2};1,\tfrac{7}{4};1\right)=\frac{3}{\pi}$ or to computing\\ 
$\int_{0}^{1}x^{-1/4}(1-x)K(x)\,dx$. A similar problem concerns the explicit evaluation of
\begin{eqnarray*} \sum_{n\geq 0}\frac{1}{(2n+1)^2}\left[\frac{1}{4^n}\binom{2n}{n}\right]^2 &=& \phantom{}_4 F_3\left(\tfrac{1}{2},\tfrac{1}{2},\tfrac{1}{2},\tfrac{1}{2};1,\tfrac{3}{2},\tfrac{3}{2};1\right)\\&\stackrel{\text{Taylor}}{=}& -\frac{2}{\pi}\int_{0}^{1} K(x^2)\log(x)\,dx.\end{eqnarray*}
We may notice that $ K(x) = \sum_{n\geq 0}P_n(2x-1) \frac{2}{2n+1}$ and
$$ -\frac{\log x}{\sqrt{x}} = \sum_{n\geq 0} P_n(2x-1) (-1)^{n} \left[8\sum_{k=0}^{n}\frac{1}{2k+1}-\frac{4}{2n+1}\right], $$
hence we may convert the original series into a combination of Euler sums:
\begin{eqnarray*} \sum_{n\geq 0}\frac{1}{(2n+1)^2}\left[\frac{1}{4^n}\binom{2n}{n}\right]^2 &=& \frac{1}{\pi}\sum_{n\geq 0}(-1)^n \frac{8\left(H_{2n+1}-\frac{1}{2}H_{2n}\right)-\frac{4}{2n+1}}{(2n+1)^2}\\
 &=& \frac{\pi^2}{8}+\frac{8}{\pi}\sum_{n\geq 1}\frac{H_{2n}(-1)^n}{(2n+1)^2}-\frac{4}{\pi}\sum_{n\geq 1}\frac{H_n(-1)^n}{(2n+1)^2},
\end{eqnarray*}
then by invoking $\sum_{n\geq 1}H_n x^n = \frac{-\log(1-z)}{1-z}$ and $\int_{0}^{1}z^k \log(z)\,dz = -\frac{1}{(k+1)^2}$ we immediately have
\begin{eqnarray*} \sum_{n\geq 1} \frac{(-1)^n H_{n}}{(2n+1)^2} &=& \int_{0}^{1}\frac{\log(1+z^2)\log(z)}{1+z^2}\,dz\\ &=& 2\,\text{Im}\,\text{Li}_3\left(\tfrac{1+i}{2}\right)-K\log(2)-\frac{\pi^3}{64}-\frac{\pi}{16}\log^2(2), \end{eqnarray*}
\begin{eqnarray*} \sum_{n\geq 1} \frac{(-1)^n H_{2n}}{(2n+1)^2} &=&\int_{0}^{1}\frac{\log(z)}{2}\left(\frac{\log(1-iz)}{1-iz}+\frac{\log(1+iz)}{1+iz}\right)\,dz\\ &=& -\text{Im}\,\text{Li}_3\left(\tfrac{1+i}{2}\right)-\frac{1}{2}K\log(2)+\frac{3\pi^3}{128}+\frac{\pi}{32}\log^2(2) \end{eqnarray*}
through the functional identities for $\text{Li}_2$ and $\text{Li}_3$. 
The terms related to $K\log(2)$ and $\pi\log^2(2)$ simplify and this leads to the unexpected identity:
\begin{eqnarray*}\label{wow1} \sum_{n\geq 0}\frac{1}{(2n+1)^2}\left[\frac{1}{4^n}\binom{2n}{n}\right]^2 & = & \phantom{}_4 F_3\left(\tfrac{1}{2},\tfrac{1}{2},\tfrac{1}{2},\tfrac{1}{2};1,\tfrac{3}{2},\tfrac{3}{2};1\right)\\ &=& \frac{3\pi^2}{8}+\frac{\log^2(2)}{2}-\frac{16}{\pi}\,\text{Im}\,\text{Li}_3\left(\tfrac{1+i}{2}\right).\end{eqnarray*}
Since $\sum_{n\geq 0}\frac{1}{(2n-1)^2}\left[\frac{1}{4^n}\binom{2n}{n}\right]^2 = \frac{4}{\pi}$ and $\sum_{n\geq 0}\frac{1}{(2n-1)^3}\left[\frac{1}{4^n}\binom{2n}{n}\right]^2 = \frac{4K-6}{\pi}$,\\ by reindexing we also get:
\begin{eqnarray*}\label{wow2} \sum_{n\geq 0}\frac{1}{(2n-1)^4}\left[\frac{1}{4^n}\binom{2n}{n}\right]^2 &=& \phantom{}_4 F_3\left(-\tfrac{1}{2},-\tfrac{1}{2},-\tfrac{1}{2},-\tfrac{1}{2};1,\tfrac{1}{2},\tfrac{1}{2};1\right) \\ &=& \frac{3\pi^2}{8}+\frac{\log^2(2)}{2}-\frac{16}{\pi}\,\text{Im}\,\text{Li}_3\left(\tfrac{1+i}{2}\right)-\frac{8}{\pi}(K-1).\end{eqnarray*}
The LHS of (\ref{wow1}) is well-approximated by the more elementary
$$ 1+\frac{1}{\pi}\sum_{n\geq 1}\frac{1}{(2n+1)^2\left(n+\frac{1}{4}\right)} = 3-\frac{\pi}{2}-\frac{4}{\pi}+\frac{4}{\pi}\log(2)$$
and similar approximated identities can be produced by replacing $\left[\frac{1}{4^n}\binom{2n}{n}\right]^2$ with $\frac{1}{\pi\left(n+\frac{1}{4}\right)}$ 
in the series representing $\left\| f(x) \right\|_{L^2(0,1)}^2$ for $f(x)=\left[x(1-x)\right]^{k/2}$ and $k$ being an odd natural number. As an alternative, 
using Theorem \ref{FLdue},
$$ \frac{16}{3\pi}=\sum_{n\geq 0}\frac{\left[\frac{1}{4^n}\binom{2n}{n}\right]^2}{(n+1)\left(\frac{1}{2}-n\right)}\approx 2+\frac{1}{\pi}\sum_{n\geq 1}\frac{1}{(n+1)\left(\frac{1}{2}-n\right)\left(n+\frac{1}{4}\right)}$$
implies $\pi\approx 4-\frac{16}{13}\log(2)$.
Many series involving squared central binomial coefficients, and eventually harmonic numbers, can be managed in a similar fashion. 
The following series appear in \cite{Campbell1},\cite{Campbell2},\cite{Parbelos}, too:
\begin{eqnarray*} \sum_{n\geq 0}\frac{H_{n+1}}{n+1}\left[\frac{1}{4^n}\binom{2n}{n}\right]^2 &=& -\frac{2}{\pi}\int_{0}^{1}K(x)\log(1-x)\,dx\\[-0.2cm]
 &=& \frac{2}{\pi}\left[2+\sum_{n\geq 1}\frac{2}{(2n+1)^2}\left(\frac{1}{n}+\frac{1}{n+1}\right)\right]\\
 &=& \frac{16}{\pi}(1-\log 2)\\
 &\approx& \frac{17}{9}+\frac{\pi}{6}-\log 4+\frac{16}{3\pi}(K+\log 2)+\frac{2}{\pi}\left(2+3\log^2 2\right),
\end{eqnarray*} 
\begin{eqnarray*} 
\!\!\!\!\!\!\!\!\!\!\!\!\!\!\!\!\!\!\!\!\!\!\!\!\!\!\!\! \sum_{n\geq 0}\frac{1}{(n+1)^2}\left[\frac{1}{4^n}\binom{2n}{n}\right]^2 &=& -\frac{2}{\pi}\int_{0}^{1}K(x)\log(x)\,dx\\
 &=& \frac{2}{\pi}\left[2+\sum_{n\geq 1}\frac{2(-1)^n}{(2n+1)^2}\left(\frac{1}{n}+\frac{1}{n+1}\right)\right]\\
 &=& -4+\frac{16}{\pi}\\
 &\approx & \frac{17}{9}-\frac{4}{\pi }-\frac{2 \pi }{9}+\frac{16 \log(2)}{3 \pi }.
\end{eqnarray*}

\begin{theorem}
We may notice that the computation of the FL-expansion of $\frac{\log(x)^{k}}{x}$ and $\frac{\log(x)^k}{\sqrt{x}}$ allows to convert 
\emph{any} hypergeometric series of the\\ $\sum_{n\geq 0}\frac{1}{(2n+1)^m}\left[\frac{1}{4^n}\binom{2n}{n}\right]^2$ or 
$\sum_{n\geq 0}\frac{1}{(n+1)^m}\left[\frac{1}{4^n}\binom{2n}{n}\right]^2$ kind into a combination of Euler sums, since the derivatives 
of the initial functions are immediately related to Stirling numbers of the first kind and central binomial coefficients.
\end{theorem}
For instance, the FL-expansion of $\frac{\log^2(x)}{\sqrt{x}}$ can be computed by noticing that
$$\frac{d^n}{dx^n}\left(\frac{\log^2 x}{\sqrt{x}}\right) = \frac{(-1)^n}{x^n\sqrt{x}}\left[A_n - B_n \log(x) + C_n \log^2(x)\right]$$
where trivially $C_n = \frac{(2n-1)!!}{2^n}$ and $B_n = \frac{2n-1}{2}B_{n-1}+2 C_{n-1}$. By letting\\ $B_n=\frac{(2n-1)!!}{2^n}\widehat{B}_n$ 
we get $\widehat{B}_n = \widehat{B}_{n-1}+\frac{4}{2n-1}$, hence 
$ B_n = \frac{(2n-1)!!}{2^n}\sum_{k=1}^{n}\frac{4}{2k-1}.$
Similarly $A_n=\frac{2n-1}{2}A_{n-1}+B_{n-1}$ leads to
\begin{eqnarray*} A_n &=& \frac{(2n-1)!!}{2^n}\cdot 8\!\sum_{1\leq k< m\leq n}\frac{1}{(2k-1)(2m-1)}\\ 
                      &=& \frac{(2n-1)!!}{2^n}\cdot 4\left[\left(\sum_{k=1}^{n}\frac{1}{2k-1}\right)^2-\sum_{k=1}^{n}\frac{1}{(2k-1)^2}\right]\end{eqnarray*}
hence by Rodrigues' formula the FL expansion of $\frac{\log^2 x}{\sqrt{x}}$ is given by:
\begin{eqnarray*} && \sum_{n\geq 0} (2n+1)P_n(2x-1)\int_{0}^{1}\frac{(1-x)^n}{n!\sqrt{x}}\left[A_n-B_n\log(x)+C_n\log^2(x)\right]\,dx=\\
%&=& \sum_{n\geq 0}(-1)^n P_n(2x-1)\left[32\left(\sum_{k=1}^{n}\frac{1}{2k-1}\right)^2+\frac{32}{2n+1}\sum_{k=1}^{n}\frac{1}{2k-1}+\frac{16}{(2n+1)^2}\right]\\
&& \sum_{n\geq 0}(-1)^n P_n(2x-1)\left[32\left(\sum_{k=0}^{n}\frac{1}{2k+1}\right)^2\!\!-\frac{32}{2n+1}\sum_{k=0}^{n}\frac{1}{2k+1}+\frac{16}{(2n+1)^2}\right].\end{eqnarray*}
This approach leads to the equality between $\sum_{n\geq 0}\frac{1}{(2n+1)^3}\left[\frac{1}{4^n}\binom{2n}{n}\right]^2 $ and 
\begin{eqnarray*} &&\frac{1}{\pi}\int_{0}^{1}K(x^2)\log^2(x)\,dx = \frac{1}{8\pi}\int_{0}^{1}K(x)\frac{\log^2(x)}{\sqrt{x}}\,dx\\
&=&\frac{8}{\pi}\sum_{n\geq 0}\frac{(-1)^n\left(\sum_{k=0}^{n}\frac{1}{2k+1}\right)^2}{(2n+1)^2}-\frac{8}{\pi}\sum_{n\geq 0}\frac{(-1)^n\sum_{k=0}^{n}\frac{1}{2k+1}}{(2n+1)^3}+\frac{4}{\pi}\sum_{n\geq 0}\frac{(-1)^n}{(2n+1)^4}.\end{eqnarray*}
The evaluation of the original series is now shown to be equivalent to the evaluation of standard and alternating Euler sums with weight $4$. We trivially have $\frac{4}{\pi}\sum_{n\geq 0}\frac{1}{(2n+1)^4}=\frac{\pi^3}{24}$ and since $\sum_{n\geq 1}H_n x^n=-\frac{\log(1-x)}{1-x}$ we also have:
\begin{eqnarray*}&&\sum_{n\geq 0}\frac{\sum_{k=0}^{n}\frac{1}{2k+1}}{(2n+1)^3} = \sum_{n\geq 0}\frac{H_{2n+1}}{(2n+1)^3}-\frac{1}{2}\sum_{n\geq 0}\frac{H_n}{(2n+1)^3}\\
&=& \sum_{n\geq 1}\frac{H_n}{n^3}-\sum_{n\geq 1}\frac{H_{2n}}{(2n)^3}-\frac{1}{2}\sum_{n\geq 0}\frac{H_n}{(2n+1)^3}\\
&=& \frac{\pi^4}{72}-\sum_{n\geq 1}\frac{H_{2n}}{(2n)^3}+\frac{1}{4}\int_{0}^{1}\frac{\log(1-x^2)\log^2(x)}{1-x^2}\,dx\\
&=& \frac{7\pi^4}{1152}+\frac{7}{8}\log(2)\zeta(3)+\int_{0}^{1}\left[\frac{\log(1-x)}{1-x}+\frac{\log(1+x)}{1+x}\right]\frac{\log^2(x)}{4x}\,dx\\
&=& -\frac{\pi^4}{1152}+\frac{7}{8}\log(2)\zeta(3)+\int_{0}^{1}\frac{\log(1+x)\log^2(x)}{4x(1+x)}\,dx\\
&=& \frac{23\pi^4}{5760}+\frac{7}{8}\log(2)\zeta(3)-\int_{0}^{1}\frac{\log(1+x)\log^2(x)}{4(1+x)}\,dx\\
&=& \frac{83\pi^4}{5760}+\frac{\pi^2}{24}\log^2(2)-\frac{\log^4(2)}{24}-\operatorname{Li}_4\left(\tfrac{1}{2}\right),\end{eqnarray*}
in agreement with De Doelder \cite{DeDoelder} and Sitaramachandra Rao \cite{Sita}. 
Once we define $\mathscr{H}_n^{(s)}=\sum_{k=0}^{n}\frac{1}{(2k+1)^s}$, in order to evaluate the series $\sum_{n\geq 0}\frac{\mathscr{H}_n^2}{(2n+1)^2}$ it is enough 
to exhibit a slight generalization of the Theorem 4.1 of Flajolet and Salvy \cite{Flajolet}. Though the residues of $\left[\psi(-s)+\gamma\right]^3$ 
they prove the following reduction formula for quadratic Euler sums:
$$ S_{1^2,q} = S_{2,q}+ q\,S_{1,q+1} - \frac{q(q+1)}{6}\zeta(q+2) + \zeta(2)\zeta(q) $$
which restricted to the case $q=2$, by symmetry, simply leads to
$$ S_{1^2,2} = 2\, S_{1,3} +\frac{13\pi^4}{360}.$$
In our case we may notice that, in the same spirit, $\sum_{n\geq 0}\frac{\mathscr{H}_n^2}{(2n+1)^2}$ equals
\begin{eqnarray*}&& \sum_{0\leq j,k \leq n}\frac{1}{(2n+1)^2(2k+1)(2j+1)}\\
&=& \sum_{0\leq j\leq n}\frac{1}{(2n+1)^2(2j+1)^2}+2\sum_{0\leq j < k\leq n}\frac{1}{(2n+1)^2 (2k+1)(2j+1)} \\
&=& \frac{1}{2}\left[\mathscr{H}_{\infty}^{(2)}\cdot \mathscr{H}_{\infty}^{(2)}+\mathscr{H}_{\infty}^{(4)}\right]+\sum_{0\leq j<k\leq n}\frac{\frac{1}{2j+1}-\frac{1}{2k+1}}{(2n+1)^2(k-j)}\\
&=& \frac{5\pi^4}{384}+\sum_{n\geq 0}\frac{1}{(2n+1)^2}\sum_{j=0}^{n}\frac{H_{n-j}-H_j}{2j+1}=\frac{5\pi^4}{384}+\mathscr{Z}(1,1,2). 
\end{eqnarray*}
For any $n\in\mathbb{N}$ we have
$$ \mathscr{H}_n^2 = -\operatorname*{Res}_{z=-(2n+1)}\frac{1}{4}\left(\gamma+2\log 2+\psi\left(\frac{z}{2}\right)\right)^2\left(\gamma+\psi(s)\right) $$
since for any $n\in\mathbb{N}$ the Laurent series of $\gamma+\psi(z)$ centered at $z=-n$ is given by $-\frac{1}{z+n}+H_n+\left(H_n^{(2)}+\zeta(2)\right)(z+n)+O((z+n)^2)$ and, 
similarly, for any even $n\in\mathbb{N}$ the Laurent series of $\frac{1}{2}\left(\gamma+2\log 2+\psi\left(\frac{z}{2}\right)\right)$ centered at $z=-n$ is given by
$$-\frac{1}{z+n}+\frac{1}{2}\left(H_{n/2}+2\log(2)\right)+\frac{1}{4}\left(H_{n/2}^{(2)}+\zeta(2)\right)(z+n)+O((z+n)^2).$$
Therefore the computation of $\sum_{n\geq 0}\frac{\mathscr{H}_n^2}{(2n+1)^2}$ can be performed by considering the residues of  
$\frac{1}{4z^2}\left(\gamma+2\log 2+\psi\left(\frac{z}{2}\right)\right)^2\left(\gamma+\psi(z)\right)$, which boils down to 
computing the residue at the origin ($\frac{\pi^2}{4}-\log^2(2)$) and the Euler sums $\sum_{n\geq 1}\frac{H_n^2}{n^2}=\frac{17\pi^4}{360},\sum_{n\geq 1}\frac{H_n^{(2)}}{n^2}=\frac{7\pi^4}{360}$, 
together with the less elementary $\sum_{n\geq 1}\frac{H_n H_{2n}}{n^2}$.\\
On the other hand, we may also consider that, by continuing $\mathscr{H}_n$ to negative integers through $\mathscr{H}_{-1}=0$ and $\mathscr{H}_{-(n+2)}=\mathscr{H}_n$,
\begin{eqnarray*} \sum_{n\in\Z}\frac{\mathscr{H}_n^2}{(2n+1)^2}&=&2\sum_{n\geq 0}\frac{\mathscr{H}_n^2}{(2n+1)^2}-2\sum_{n\geq 0}\frac{\mathscr{H}_n}{(2n+1)^3}+\sum_{n\geq 0}\frac{1}{(2n+1)^4}\\
&=&\frac{1}{4}\sum_{n\in\Z}\frac{\left(\psi\left(n+\frac{3}{2}\right)-\psi\left(\frac{1}{2}\right)\right)^2}{(2n+1)^2}\\
&=&-\frac{1}{4}\sum_{m\geq 0}\operatorname*{Res}_{z=-(2m+1)/2} \pi\cot(\pi z)\left(\frac{\psi\left(z+\frac{3}{2}\right)-\psi\left(\frac{1}{2}\right)}{2z+1}\right)^2\\
&=& \frac{\pi^2}{4}\log^2(2)+\frac{\pi^4}{96}.
\end{eqnarray*}
This identity has many direct consequences:
$$\mathscr{Z}(1,1,2) = \frac{\pi^4}{720}+\frac{\pi^2}{6}\log^2(2)-\frac{\log^4(2)}{24}-\text{Li}_4\left(\tfrac{1}{2}\right),$$
\begin{eqnarray*}\mu_1 &=& \frac{1}{2}\int_{0}^{1}\frac{\log^2(z)\log(1+z)}{z(1+z)}\,dz= \sum_{n\geq 1}\frac{(-1)^{n+1} H_n}{n^3}\\&=& \frac{11\pi^4}{360}+\frac{\pi^2}{12}\log^2(2)-\frac{\log^4(2)}{12}-\frac{7}{4}\log(2)\zeta(3)-2\,\text{Li}_4\left(\tfrac{1}{2}\right),\end{eqnarray*}
\begin{eqnarray*}\sum_{n\geq 0}\frac{H_{n+1}^2+H_{n+1}^{(2)}}{n+1}\left[\frac{1}{4^n}\binom{2n}{n}\right]^2&=&\frac{2}{\pi}\int_{0}^{1}K(x^2)\log^2(1-x)\,dx\\&=&\frac{96}{\pi}-\frac{8\pi}{3}-\frac{128}{\pi}\log(2)+\frac{64}{\pi}\log^2(2)\end{eqnarray*}
and it also allows the explicit evaluation of 
\begin{eqnarray*} \sum_{n\geq 1}\frac{H_n H_{2n}}{n^2}&=&2\!\!\sum_{a,b,n\geq 1}\frac{1}{ab(a+n)(b+2n)}\\
&=&2\iint_{(0,1)^2}\frac{w\log(1-w)\log(1-z)}{1-w^2 z}\,dw\,dz\end{eqnarray*}
in terms of $\pi,\log(2),\zeta(3)$ and $\text{Li}_4\left(\tfrac{1}{2}\right)$. We are now ready to tackle the alternating Euler sums with weight $4$ involving the powers of $\mathscr{H}_n$. By mimicking the approach used to 
evaluate $\sum_{n\geq 0}\frac{\mathscr{H}_n^2}{(2n+1)^2}$ we have that
\begin{eqnarray*} &&\sum_{n\in\mathbb{Z}}\frac{(-1)^n \mathscr{H}_n^2}{(2n+1)^2} =  2 \sum_{n\geq 0}\frac{(-1)^n \mathscr{H}_n}{(2n+1)^3}-\sum_{n\geq 0}\frac{(-1)^n}{(2n+1)^4}\\
&=& -\frac{1}{4}\sum_{m\geq 0}\operatorname*{Res}_{z=-(2m+1)/2}\frac{\pi}{\sin(\pi z)}\left(\frac{\psi\left(z+\frac{3}{2}\right)-\psi\left(\frac{1}{2}\right)}{2z+1}\right)^2\\
&=& \frac{\pi^3}{24}\log(2)-\frac{1}{4}\sum_{m\geq 1}\left[-\frac{\pi(-1)^m}{2m^3}+\frac{\pi(-1)^m}{2m^2}\left(H_{m-1}^{(2)}+2\log 2\right)\right]\\
&=& \frac{\pi^3}{16}\log(2)-\frac{3\pi}{32}\zeta(3)-\frac{\pi}{8}\int_{0}^{1}\frac{\text{Li}_2(-x)\log(x)}{1+x}\,dx
\end{eqnarray*}
where the last integral depends on $\mu_1$: it equals $\frac{13\pi^4}{288}-\frac{\log^4(2)}{6}+\frac{\pi^2}{6}\log^2(2)-4\,\text{Li}_4\left(\tfrac{1}{2}\right)-\frac{7}{2}\log(2)\zeta(3)$. Similarly
\begin{eqnarray*} &&\sum_{n\in\mathbb{Z}}\frac{(-1)^n \mathscr{H}_n^3}{2n+1}\\&=&2 \sum_{n\geq 0}\frac{(-1)^n \mathscr{H}_n^3}{2n+1}-3 \sum_{n\geq 0}\frac{(-1)^n \mathscr{H}_n^2}{(2n+1)^2}+3\sum_{n\geq 0}\frac{(-1)^n\mathscr{H}_n}{(2n+1)^3}-\sum_{n\geq 0}\frac{(-1)^n}{(2n+1)^4}\\
&=& -\sum_{m\geq 0}\operatorname*{Res}_{z=-(2m+1)/2}\frac{\pi}{\sin(\pi z)(2z+1)}\left(\frac{\psi\left(z+\frac{3}{2}\right)-\psi\left(\frac{1}{2}\right)}{2}\right)^3\\
&=& \frac{\pi}{2}\log^3(2)-\frac{\pi}{8}\sum_{m\geq 1}(-1)^m\left[-\frac{6}{m}\left(\frac{H_{m-1}}{2}+\log 2\right)^2+\frac{3}{m}\left(\frac{H_{m-1}^{(2)}}{2}+\frac{\pi^2}{12}\right)\right]\\
&& -\frac{\pi}{8}\sum_{m\geq 1}(-1)^m\left[\frac{3}{m^2}\left(\frac{H_{m-1}}{2}+\log 2\right)-\frac{1}{2m^3}-\frac{\pi^2}{4m}\right]\\
%&=& \frac{\pi}{2}\log^3(2)-\frac{\pi}{8}\sum_{m\geq 1}(-1)^m\left[-\frac{5}{m^3}+\frac{9H_m}{2 m^2}-\frac{3(H_m^2-H_m^{(2)})}{2 m}+\frac{9\log(2)}{m^2}-\frac{6 H_m \log(2)}{m}-\frac{6 \log^2(2)}{m}\right]\\
&=& \frac{\pi}{8}\log^3(2)+\frac{\pi^3}{32}\log(2)-\frac{15\pi}{128}\zeta(3)+\frac{3\pi}{16}\sum_{m\geq 1}\frac{(-1)^m}{m}\left(H_m^2-H_m^{(2)}\right)\\
&=& \frac{\pi}{8}\log^3(2)+\frac{\pi^3}{32}\log(2)-\frac{15\pi}{128}\zeta(3)+\frac{3\pi}{16}\int_{0}^{1}\frac{\log^2(1+x)}{x(1+x)}\,dx\\
&=& \frac{\pi}{16}\log^3(2) +\frac{\pi^3}{32}\log(2)-\frac{9\pi}{128}\zeta(3),
\end{eqnarray*}
such that 
\begin{eqnarray*}\label{wow3}&&\sum_{n\geq 0}\frac{1}{(2n+1)^3}\left[\frac{1}{4^n}\binom{2n}{n}\right]^2 =\phantom{}_5 F_4\left(\tfrac{1}{2},\tfrac{1}{2},\tfrac{1}{2},\tfrac{1}{2},\tfrac{1}{2};1,\tfrac{3}{2},\tfrac{3}{2},\tfrac{3}{2};1\right)\\&=& \frac{16}{3\pi}\sum_{n\geq 0}\frac{(-1)^n \mathscr{H}_n^3}{2n+1}-\frac{\log^3(2)}{6}-\frac{\pi^2}{12}\log(2)+\frac{3}{16}\zeta(3)+\frac{4}{3\pi}\beta(4).\end{eqnarray*}
The previous equation gives an explicit relation between a cubic, alternating Euler sum with weight four and a quadratic, alternating sum with the same weight. 
These series are further generalizations of the quadratic, alternating Euler sums with weight four computed by Xu, Yang and Zhang in \cite{XYZ}. 
The combination $2\sum_{n\in\Z}\frac{(-1)^n \mathscr{H}_n^3}{2n+1}-3\sum_{n\in\Z}\frac{(-1)^n \mathscr{H}_n^2}{(2n+1)^2}$ has been computed by Zheng in \cite{Zeng} ($(4.5a),(4.5b)$), by exploiting the Dougall-Dixon theorem 
and by comparing the coefficients of suitable infinite products. In the FL-context we may notice that
\begin{equation}\label{FLstorta} 2\,K(x)-\frac{\log(1-x)}{\sqrt{1-x}} \stackrel{L^2(0,1)}{=} 8\sum_{n\geq 0}\mathscr{H}_n P_n(2x-1) \end{equation}
immediately leads to
\begin{eqnarray*}\label{FLs2}\notag 64 \sum_{n\geq 0}\frac{(-1)^n \mathscr{H}_n^2}{2n+1}&=&\frac{\pi^3}{3}+4\pi\log^2(2)-4\int_{0}^{1}K(x)\frac{\log x}{\sqrt{x}}\,dx\\
&=& \frac{\pi^3}{3}+4\pi\log^2(2)+8\pi \sum_{n\geq 0}\frac{1}{(2n+1)^2}\left[\frac{1}{4^n}\binom{2n}{n}\right]^2  \\
&=& \frac{10\pi^3}{3}+8\pi\log^2(2)-128\,\text{Im}\,\text{Li}_3\left(\tfrac{1+i}{2}\right).
\end{eqnarray*}
By exploiting the expansion of $\frac{\log(x)}{\sqrt{1-x}}$ it is possible to compute $\sum_{n\geq 0}\frac{(-1)^n \mathscr{H}_n^2}{(2n+1)^3}$\\ in a similar way.
Proving the following identities is a simpler task: 
\begin{eqnarray*}&&\sum_{n\geq 0}\frac{1}{(n+1)^3}\left[\frac{1}{4^n}\binom{2n}{n}\right]^2=\phantom{}_4 F_3\left(\tfrac{1}{2},\tfrac{1}{2},1,1;2,2,2;1\right)\\&=&\frac{16}{\pi}\left[3-2K-\pi+\pi\log 2\right],\end{eqnarray*}
\begin{eqnarray*}&&\sum_{n\geq 0}\frac{1}{(n+1)^4}\left[\frac{1}{4^n}\binom{2n}{n}\right]^2=\phantom{}_5 F_4\left(\tfrac{1}{2},\tfrac{1}{2},1,1,1;2,2,2,2;1\right)\\&=&\frac{2}{\pi}\left[64-64K-24\pi+3\pi^3-128\,\text{Im}\,\text{Li}_3\left(\tfrac{1+i}{2}\right)+32\,\pi\log 2-12\,\pi\log^2 2\right].\end{eqnarray*}
They can be deduced in a straightforward way from the FL-expansions of $\frac{d}{dx}E(x)$ (computed in our previous work \cite{Interplay}), $\log(x)$ and $\log^2(x)$. They also lead to
$$ \sum_{n\geq 0}\frac{1}{(2n-1)^2}\left[\frac{1}{4^n}\binom{2n}{n}\right]^2 = \phantom{}_2 F_1\left(-\tfrac{1}{2},-\tfrac{1}{2};1;1\right)=\frac{4}{\pi} $$
$$ \sum_{n\geq 0}\frac{1}{(2n-1)^2 (2n-3)^2}\left[\frac{1}{4^n}\binom{2n}{n}\right]^2 = \tfrac{1}{9}\cdot \phantom{}_2 F_1\left(-\tfrac{3}{2},-\tfrac{3}{2};1;1\right)=\frac{32}{27\pi} $$
consistent with the representations of $\frac{1}{\pi}$ obtained in \cite{Interplay} by computing the moments of generalized elliptic integrals.
Similar manipulations might be useful in improving the state of the art about the irrationality measure of $\pi$, since it is much simpler to estimate the irrationality measure of $\log(2)$ 
through Beuker-like integrals and Viola's method, as done by Viola himself and Marcovecchio in \cite{Viola}. Hypergeometric functions extending the ``usual'' 
kernels $x^a(1-x)^b$ and $P_n(x)$ have been crucial in the work of Zudilin \cite{Zudilin} about the irrationality measure of $\pi^2$, too. An application of the Cauchy-Schwarz inequality leads to
\begin{eqnarray*}2-\frac{16}{3\pi}&=&\sum_{n\geq 1}\frac{\left[\frac{1}{4^n}\binom{2n}{n}\right]^2}{(n+1)\left(n-\frac{1}{2}\right)}\\&\stackrel{\text{CS}}{\leq}&\sqrt{\sum_{n\geq 1}\frac{1}{(n+1)\left(n-\frac{1}{2}\right)\left(n+\frac{1}{4}\right)}\sum_{n\geq 1}\frac{\left(n+\frac{1}{4}\right)\left[\frac{1}{4^n}\binom{2n}{n}\right]^4}{(n+1)\left(n-\frac{1}{2}\right)}}\\
 &=& \sqrt{\frac{8}{9}\left(7-\pi-4\log 2\right)\left(\frac{1}{2}-\frac{4}{\pi^2}\right)}
\end{eqnarray*}
where the series involving $\left[\frac{1}{4^n}\binom{2n}{n}\right]^4$ has been computed by considering the FL-expansions of $\sqrt{x(1-x)}$ and $\frac{1}{\sqrt{x(1-x)}}$.
This already leads to $\pi\approx 3.141609$ and it can be generalized by replacing the $(n+1)\left(\frac{1}{2}-n\right)$ term with $(n+1)_k\left(\frac{1}{2}-n\right)_k$ for $k>1$.
Approximations of zero through linear combinations with rational coefficients of $1,\pi,\pi^3,\log(2)$ and $\zeta(3)$ can be produced by considering the more accurate
estimation $\left[\frac{1}{4^n}\binom{2n}{n}\right]^2 \approx \frac{1}{\pi\left(n+\frac{1}{4}\right)}-\frac{1}{32\pi\left(n+\frac{1}{4}\right)^3}$ as $n\gg 1$.\\
We may also consider that by creative telescoping, the following identity holds for any $n\in\N$:
\begin{eqnarray}\label{CT1}\notag &&\sum_{m\geq 1}\frac{\Gamma\left(m-\frac{1}{4}\right)^2}{\Gamma\left(-\frac{1}{4}\right)^2 (n+1)_m \left(\frac{1}{2}-n\right)_m}\\&=& -\frac{1}{(4n+1)^2}+(-1)^n \frac{4^n\sqrt{\pi}}{(4n+1)^2\,\Gamma\left(\frac{3}{4}\right)^2\binom{2n}{n}}.\end{eqnarray}
If we multiply both sides by $(-1)^n (4n+1)\left[\frac{1}{4^n}\binom{2n}{n}\right]^3$ and sum over $n\geq 0$ we get:
\begin{eqnarray}\label{CT2}\notag&&\sum_{m\geq 1}\frac{\Gamma\left(m-\frac{1}{4}\right)^2}{\Gamma\left(-\frac{1}{4}\right)^2}\sum_{n\geq 0}\frac{(-1)^n (4n+1)\left[\frac{1}{4^n}\binom{2n}{n}\right]^3}{(n+1)_m\left(\frac{1}{2}-n\right)_m}\\&=& \frac{\Gamma\left(\frac{1}{4}\right)^6}{32\pi^{7/2}}-\phantom{}_4 F_3\left(\tfrac{1}{4},\tfrac{1}{2},\tfrac{1}{2},\tfrac{1}{2};1,1,\tfrac{5}{4};-1\right)\end{eqnarray}
and by exploiting the consequence of Parseval's identity provided by (\ref{nove}) we get:
\begin{eqnarray}\label{CT3}\notag\frac{2}{\pi}\sum_{m\geq 1}\frac{\Gamma\left(m-\frac{1}{4}\right)^2 16^m m!^2}{\Gamma\left(-\frac{1}{4}\right)^2 (2m)!^2}&=& \frac{1}{2\pi}\phantom{}_3 F_2\left(\tfrac{3}{4},\tfrac{3}{4},1;\tfrac{3}{2},\tfrac{3}{2};1\right)\\&=&\frac{\Gamma\left(\frac{1}{4}\right)^6}{32\pi^{7/2}}-\phantom{}_4 F_3\left(\tfrac{1}{4},\tfrac{1}{2},\tfrac{1}{2},\tfrac{1}{2};1,1,\tfrac{5}{4};-1\right)\end{eqnarray}
reducing the problem of evaluating a $\phantom{}_4 F_3\left(\ldots,-1\right)$ with quarter-integer parameters to the problem of evaluating a $\phantom{}_3 F_2(\ldots,1)$ with 
quarter-integer parameters. On its turn the last problem is equivalent to finding the $L^2$ norm over $S^1$ for $\phantom{}_2 F_1\left(\tfrac{3}{4},1;\tfrac{3}{2};z\right)=\int_{0}^{1}\frac{dt}{(1-z+z t^2)^{3/4}}$. In the opposite direction the identity above can be seen as an acceleration formula, since the main term of the series defining
$\phantom{}_4 F_3(\ldots,-1)$ has an alternating sign and a faster decay to zero with respect to the main term of the series defining $\phantom{}_3 F_2(\ldots ,1)$.\\
If we multiply both sides of (\ref{CT1}) by $\left[\frac{1}{4^n}\binom{2n}{n}\right]^2$ and sum over $n\geq 0$, by Theorem (\ref{FLdue}) we get:
\begin{eqnarray*}\label{CT4}\sum_{m\geq 1}\frac{\Gamma\left(m-\frac{1}{4}\right)^2}{\Gamma\left(-\frac{1}{4}\right)^2}\sum_{n\geq 0}\frac{\left[\frac{1}{4^n}\binom{2n}{n}\right]^2}{(n+1)_m\left(\frac{1}{2}-n\right)_m}&=& \frac{1}{3\pi}\phantom{}_4 F_3\left(\tfrac{3}{4},\tfrac{3}{4},1,1;\tfrac{5}{4},\tfrac{3}{2},\tfrac{7}{4};1\right)\\
&=&\frac{1}{8\pi}\int_{0}^{1}K(x)x^{-3/4}\log(x)\,dx \\ &-& \frac{\sqrt{\pi}}{4\,\Gamma\left(\frac{3}{4}\right)^2}\int_{0}^{1}\frac{\log(x)\,dx}{\sqrt{x(1+x^2)}}\\
&=& \frac{\sqrt{\pi}}{\Gamma\left(\frac{3}{4}\right)^2}\phantom{}_3 F_2\left(\tfrac{1}{4},\tfrac{1}{4},\tfrac{1}{2};\tfrac{5}{4},\tfrac{5}{4};-1\right)\\ &-& \phantom{}_4 F_3\left(\tfrac{1}{4},\tfrac{1}{4},\tfrac{1}{2},\tfrac{1}{2};1,\tfrac{5}{4},\tfrac{5}{4};1\right).\end{eqnarray*}
Similarly, if we multiply both sides of (\ref{CT1}) by $\left[\frac{1}{4^n}\binom{2n}{n}\right]^3$ and sum over $n\geq 0$, by Theorem (\ref{FLdue}) we get:
\begin{eqnarray*}\label{CT5}\sum_{m\geq 1}\frac{\Gamma\left(m-\frac{1}{4}\right)^2}{\Gamma\left(-\frac{1}{4}\right)^2}\sum_{n\geq 0}\frac{\left[\frac{1}{4^n}\binom{2n}{n}\right]^3}{(n+1)_m\left(\frac{1}{2}-n\right)_m}&=& \frac{\Gamma\left(\frac{1}{4}\right)^4}{48\pi^3}\phantom{}_4 F_3\left(\tfrac{3}{4},\tfrac{3}{4},1,\tfrac{5}{4};\tfrac{3}{2},\tfrac{3}{2},\tfrac{7}{4};1\right)\\
&=&\frac{1}{\Gamma\left(\frac{3}{4}\right)^2}\phantom{}_4 F_3\left(\tfrac{1}{4},\tfrac{1}{4},\tfrac{1}{2},\tfrac{1}{2};1,\tfrac{5}{4},\tfrac{5}{4};-1\right)\\  &-& \phantom{}_5 F_4\left(\tfrac{1}{4},\tfrac{1}{4},\tfrac{1}{2},\tfrac{1}{2},\tfrac{1}{2};1,1,\tfrac{5}{4},\tfrac{5}{4};1\right).\end{eqnarray*}

\section{Twisted hypergeometric series}
When manipulating hypergeometric series through the binomial or Euler's transforms, generalized harmonic numbers $H_n^{(s)}$ may appear 
among the factors of the new main term. For instance, some solutions to the Basel problem $\zeta(2)=\frac{\pi^2}{6}$ exploit the identity
$\zeta(2) = 2\sum_{n\geq 1}\frac{H_n}{n 2^n}$. In this context it is pretty natural to wonder if a technique for the explicit evaluation
of some $\phantom{}_p F_q(\ldots;\pm 1)$ allows the evaluation of \emph{twisted} series too, i.e. series of the form $\sum_{n\geq 0} H_n^{(s)} a_n$
where $\sum_{n\geq 0} a_n = \phantom{}_p F_q(\ldots;\pm 1)$. Some cases involving central binomial coefficients and $H_n$ or $H_{2n}$ have 
already been investigated by Campbell through integral transforms (see \cite{Campbell2}); this section is devoted to showing how to deal with such \emph{twisted} 
series through the FL machinery.\\
Rodrigues' and Bonnet's formulas lead to the FL-expansions
$$ \log^2(1-x) = 2 + \sum_{n\geq 1}P_n(2x-1)\left[2\left(\frac{1}{n}+\frac{1}{n+1}\right)^2+4\left(\frac{1}{n}+\frac{1}{n+1}\right)H_{n-1}\right], $$
$$ \log(x)\log(1-x) = (2-\zeta(2))-\sum_{n\geq 1}P_{2n}(2x-1)\frac{4n+1}{2n^2(2n+1)^2}, $$
$$ \frac{-\log(x)}{1-x} = \zeta(2) - \sum_{n\geq 1}P_n(2x-1)(2n+1)\left[2\sum_{k=1}^{n}\frac{(-1)^{k+1}}{k^2}-\zeta(2)\right],$$
while the generating function for Legendre polynomials gives
\begin{equation}\label{harmonicFL} \frac{1}{2\sqrt{x}} \log\left(\frac{1+\sqrt{x}}{4\sqrt{x}}\right)\stackrel{L^2(0,1)}{=}\sum_{n\geq 0}(-1)^n H_n P_n(2x-1). \end{equation}
On the other hand we have
$$ \int_{0}^{1}x^n \log(x)\log(1-x)\,dx = \frac{2}{(n+1)^3}+\frac{H_n}{(n+1)^2}-\frac{\zeta(2)-H_n^{(2)}}{n+1},$$
$$ \int_{0}^{1} x^n \log^2(1-x)\,dx = \frac{H_{n+1}^2 + H_{n+1}^{(2)}}{n+1}, \qquad  \int_{0}^{1} x^n \frac{-\log x}{1-x}\,dx = \zeta(2)-H_n^{(2)}, $$
hence these identities, through the FL machinery, allow an explicit evaluation of hypergeometric series whose main term depends on $\left[\frac{1}{4^n}\binom{2n}{n}\right]$
or $\left[\frac{1}{4^n}\binom{2n}{n}\right]^2$ and $H_n, H_n^2, H_n^{(2)}$ (hypergeometric series with a \emph{harmonic twist}, according to the 
terminology introduced by Campbell in \cite{Parbelos}). For instance

\begin{eqnarray*} \int_{0}^{1}K(x)\frac{-\log x}{1-x}\,dx &=& \frac{\pi}{2}\sum_{n\geq 0}\left(\zeta(2)-H_n^{(2)}\right)\left[\frac{1}{4^n}\binom{2n}{n}\right]^2 \\
&=& 2\,\zeta(2)+\sum_{n\geq 1}\frac{4}{2n+1}\sum_{k>n}\frac{(-1)^{k+1}}{k^2}\\
&=& 4\int_{0}^{1}\frac{\text{Li}_2(-x^2)-\text{Li}_2(-1)}{1-x^2}\,dx \\
&=& 2\int_{0}^{1}\frac{\text{Li}_2(-x)-\text{Li}_2(-1)}{(1-x)\sqrt{x}}\,dx \\
&=& 4 \int_{0}^{1}\left[\log(1+x)-\log(1-x)\right]\log(1+x^2)\frac{dx}{x}\\
&=& 4\pi K-7\zeta(3)
\end{eqnarray*}
by the functional identities for $\text{Li}_2$ and $\text{Li}_3$. Remarkably, there is a standard hypergeometric series which is closely related 
to $\sum_{n\geq 0}\left(\zeta(2)-H_n^{(2)}\right)\left[\frac{1}{4^n}\binom{2n}{n}\right]^2$:
$$ \sum_{n\geq 0}\frac{16^n}{(2n+1)^3 \binom{2n}{n}^2} = \phantom{}_4 F_3\left(\tfrac{1}{2},1,1,1;\tfrac{3}{2},\tfrac{3}{2},\tfrac{3}{2};1\right) = -\pi K+\tfrac{7}{2}\zeta(3).$$
This identity can be proved through the same technique outlined in \cite{Ditrani}: by the Maclaurin series of the $\arcsin$ function and the moments 
$ \int_{0}^{\pi/2}\left(\sin\theta\right)^{2n+1}\,d\theta=\frac{1}{2n+1}\left[\frac{1}{4^n}\binom{2n}{n}\right]^{-1}$ we have:
\begin{eqnarray*}\sum_{n\geq 0}\frac{16^n}{(2n+1)^3 \binom{2n}{n}^2}&=&\int_{0}^{\pi/2} \theta \log\tan\frac{\theta}{2}-\text{Im}\,\text{Li}_2(-e^{i\theta})+\text{Im}\,\text{Li}_2(e^{i\theta}) \,d\theta\\&\stackrel{\text{IBP}}{=}&-\pi K+\tfrac{7}{2}\zeta(3). \end{eqnarray*}
More simply,
\begin{eqnarray*} -\int_{0}^{1}K(x)\log(1-x)\,dx &\stackrel{\text{Taylor}}{=}& \frac{\pi}{2}\sum_{n\geq 0}\frac{H_{n+1}}{n+1}\left[\frac{1}{4^n}\binom{2n}{n}\right]^2 \\
&\stackrel{\text{FL}}{=}& 2+\sum_{n\geq 1}\left(\frac{1}{n}+\frac{1}{n+1}\right)\frac{2}{(2n+1)^2}\\
&=& 8(1-\log 2),
\end{eqnarray*}
\begin{eqnarray*}\!\!\! -\int_{0}^{1}K(x)\log(1-\sqrt{x})\,dx &\stackrel{\text{Taylor}}{=}& \frac{\pi}{2}\sum_{n\geq 0}\frac{H_{2n+2}}{n+1}\left[\frac{1}{4^n}\binom{2n}{n}\right]^2 \\
&\stackrel{\text{FL}}{=}& 3-\sum_{n\geq 1}\frac{2}{(2n+1)^2}\left[\frac{(-1)^n}{2n(n+1)}-\frac{2n+1}{n(n+1)}\right]\\
&=& 4+4K-6\log 2,
\end{eqnarray*}
\begin{eqnarray*} \sum_{n\geq 0}\frac{H_{2n}}{(n+1)(1-2n)}\left[\frac{1}{4^n}\binom{2n}{n}\right]^2 &\stackrel{\text{FL}}{=}& \frac{8}{\pi}\int_{0}^{1}\frac{\sqrt{x(1-x)}}{2\sqrt{x}}\log\left(\frac{1+\sqrt{x}}{4\sqrt{x}}\right)\,dx \\
&=& \frac{4}{\pi}\int_{0}^{1}\log\left(\frac{1+\sqrt{x}}{4\sqrt{x}}\right)\sqrt{1-x}\,dx \\
 \left[\text{FL of } \sqrt{1-x},\;\log\left(1+\tfrac{1}{\sqrt{x}}\right)\right]\; &=& \sum_{n\geq 0}\frac{8(-1)^n}{\pi(1-4n^2)(2n^2+5n+3)}-\frac{16\log 2}{3\pi} \\
&=& \frac{2}{3}+\frac{8}{3\pi}\left(1-\log 2\right).
\end{eqnarray*}
Of course this approach can be combined with classical approaches based on integral transforms. Three instances are provided by
the series\\ $\sum_{n\geq 1}\frac{H_n}{n}\left[\frac{1}{4^n}\binom{2n}{n}\right]$, $\sum_{n\geq 1}\frac{H_n}{n}\left[\frac{1}{4^n}\binom{2n}{n}\right]^2$ and $\sum_{n\geq 1}\frac{1}{n^2}\left[\frac{1}{4^n}\binom{2n}{n}\right]^2$:
$$\sum_{n\geq 1}\frac{H_n}{n}\left[\frac{1}{4^n}\binom{2n}{n}\right] = \zeta(2)-2\log^2(2)+\frac{1}{\pi}\int_{0}^{\pi/2}\log^2(1-\sin^2\theta)\,d\theta = 2\,\zeta(2),$$
\begin{eqnarray*}
 &&\sum_{n\geq 1}\frac{1}{n^2}\left[\frac{1}{4^n}\binom{2n}{n}\right]^2 = \phantom{}_5 F_4\left(1,1,1,\tfrac{3}{2},\tfrac{3}{2};2,2,2,2;1\right)\\[-0.2cm]
 &\stackrel{\text{Taylor}}{=}& -\frac{2}{\pi}\int_{0}^{\pi/2}\left(K(x)-\frac{\pi}{2}\right)\log(x)\frac{dx}{x}\\
 &=& \sum_{n\geq 0}\frac{(2n+1)(2(n+1)-1)}{4(n+1)^4}\left[\frac{1}{4^n}\binom{2n}{n}\right]^2\\
 &\stackrel{\text{Taylor}}{=}& \frac{1}{\pi}\int_{0}^{1}\frac{d}{dx}E(x)\cdot\left(4\log x+\log^2 x\right)\,dx\\
 &\stackrel{\text{FL}}{=}& 1-\frac{2}{\pi}+\frac{1}{\pi}\sum_{n\geq 1}(-1)^n\left[\frac{2 (2n+1)(1-2 n^2)}{n^2 (1+n)^2}+4\frac{2n+1}{n(n+1)}H_{n-1}\right]\\ &\cdot& \left[\int_{0}^{1}\frac{2x^{2n+2}}{1+x^2}\,dx-\frac{1}{2n+1}\right]\\
 &=& 1+\frac{-28+\frac{4\pi^2}{3}+32\log(2)-16\log^2(2)}{\pi}+\frac{1}{\pi}\int_{0}^{1}\frac{f(x)}{1+x^2}\,dx\\
\end{eqnarray*}
where $ f(x)= 20x^2+(24-8x^2)\log(1-x^2)+4(x^2+1)\text{Li}_2(x^2)+4(x^2-1)\log^2(1-x^2).$
In particular $\int_{0}^{1}\frac{f(x)}{1+x^2}$ is given by a linear combination with rational coefficients of $1,\pi,\log(2),K,\pi\log(2),\pi^2,\log^2(2)$ and 
$$ \int_{0}^{1}\frac{\log^2(1-x^2)}{1+x^2}\,dx =\frac{5\pi^3}{32}-4\,\text{Im}\,\text{Li}_3\left(\tfrac{1+i}{2}\right)-4K\log(2)+\frac{3\pi}{8}\log^2(2),$$
which is associated to Euler sums with weight $3$ and allows to state that\\ $\phantom{}_5 F_4\left(1,1,1,\tfrac{3}{2},\tfrac{3}{2};2,2,2,2;1\right)$ has an explicit form in terms of the mentioned constants.
The chain of equalities above is an application of the \emph{transformation formula} introduced in our previous work \cite{Interplay}, 
allowing to convert $\int_{0}^{1}\left(K(x)-\frac{\pi}{2}\right)\frac{f(x)}{x}\,dx$ into $\int_{0}^{1}\frac{d}{dx}E(x)\cdot \widetilde{f}(x)\,dx$.
The application of such principle to the series\\ $\sum_{n\geq 1}\frac{H_n}{n}\left[\frac{1}{4^n}\binom{2n}{n}\right]$ gives the following alternative derivation:
\begin{eqnarray*} \sum_{n\geq 1}\frac{H_n}{n}\left[\frac{1}{4^n}\binom{2n}{n}\right] &=& \sum_{n\geq 0}\frac{H_{n+1}}{n+1}\left[\frac{1}{4^n}\binom{2n}{n}\right]\frac{2n+1}{2n+2}\\
&\stackrel{\text{Taylor}}{=}& -\int_{0}^{1} \frac{1-\sqrt{1-x}}{x\sqrt{1-x}}\log(1-x)\,dx \\
&\stackrel{x\mapsto 1-x}{=}& -\int_{0}^{1} \frac{1-\sqrt{x}}{(1-x)\sqrt{x}}\log(x)\,dx \\
&\stackrel{x\mapsto x^2}{=}& -4\int_{0}^{1} \frac{1-x}{1-x^2}\log(x)\,dx\\
&=& 4\int_{0}^{1}\frac{-\log(x)}{1+x}\,dx = 4\sum_{n\geq 1}\frac{(-1)^{n+1}}{n^2} = 2\,\zeta(2).
\end{eqnarray*}
and the same procedure tackles
\begin{eqnarray*} &&\sum_{n\geq 1}\frac{H_n}{n 2^n}\left[\frac{1}{4^n}\binom{2n}{n}\right] \stackrel{\text{FL}}{=} \int_{0}^{1}\frac{-\log(x)}{\sqrt{1+x}\left(\sqrt{2}+\sqrt{1+x}\right)}\,dx\\ 
&=& \frac{1}{\sqrt{2}}\int_{0}^{1}\left(-\log x\right)\left(\frac{1}{\sqrt{1+x}}-\frac{1}{\sqrt{2}+\sqrt{1+x}}\right)\,dx \\
&=& \frac{1}{\sqrt{2}}\sum_{n\geq 0}\frac{(-1)^n}{(n+1)^2}\left[\frac{1}{4^n}\binom{2n}{n}\right]+\frac{1}{\sqrt{2}}\int_{0}^{1}\frac{\log(x)}{\sqrt{2}+\sqrt{1+x}}\,dx 
\end{eqnarray*}
too, where the first term of the RHS is elementary and the second term depends on $\text{Li}_2(3-2\sqrt{2})$ by the substitution $x\mapsto 1-2\left(\frac{2t}{1+t^2}\right)^2$.
The series with alternating signs $\sum_{n\geq 0}(-1)^n H_{2n}\left[\frac{1}{4^n}\binom{2n}{n}\right]$ equals $\frac{1}{\sqrt{2}}\log\left(\frac{1+\sqrt{2}}{4}\right)$ 
by evaluating $\sum_{n\geq 0} H_n P_n(2x-1)$ at $x=\frac{1}{2}$. By directly exploiting the FL-expansion of $\frac{\log(1-x)}{x}$ we have:
\begin{eqnarray*}
 &&\sum_{n\geq 1}\frac{H_n}{n}\left[\frac{1}{4^n}\binom{2n}{n}\right]^2 \stackrel{\text{Taylor}}{=} -\frac{2}{\pi}\int_{0}^{\pi/2}\left(K(x)-\frac{\pi}{2}\right)\log(1-x)\frac{dx}{x}\\
 &\stackrel{\text{FL}}{=}&\frac{2}{\pi}\left[\zeta(2)\left(2-\frac{\pi}{2}\right)+\sum_{n\geq 1}\frac{2(-1)^{n+1}}{2n+1}\left(-\zeta(2)+2\sum_{k=1}^{n}\frac{(-1)^{k+1}}{k^2}\right)\right]\\
 &=& \frac{8}{\pi}\sum_{n\geq 1}\frac{(-1)^{n+1}}{2n+1}\sum_{k=1}^{n}\frac{(-1)^{k+1}}{k^2}\\
 &=& \frac{8}{\pi}\int_{0}^{1}\frac{\text{Li}_2(x^2)}{1+x^2}\,dx\\
 &=& -\frac{5\pi ^2}{3}+\frac{64}{\pi}\, \text{Im}\,\text{Li}_3\left(\tfrac{1+i}{2}\right)+\frac{32}{\pi}K\log(2)-2\log^2(2).
\end{eqnarray*}
About the approximated series:
$$ \frac{1}{4}+\sum_{n\geq 2}\frac{H_n}{\pi n\left(n+\frac{1}{4}\right)} = \frac{1}{4}+6\log(2)+\frac{18}{\pi}\log^2(2)-\frac{\pi}{2}-\frac{4}{5\pi}-\frac{16}{\pi}K.$$
Now we may tackle two similar sums, \emph{twisted} by factors $H_n^2, H_n^{(2)}$.
\begin{eqnarray*} &&\int_{0}^{1}K(x)\log(x)\log(1-x)\,dx\\ &\stackrel{\text{Taylor}}{=}& \frac{\pi}{2}\sum_{n\geq 0}\left(\frac{H_{n+1}}{(n+1)^2}-\frac{\zeta(2)-H_{n+1}^{(2)}}{n+1}\right)\left[\frac{1}{4^n}\binom{2n}{n}\right]^2\\
&\stackrel{\text{FL}}{=}& (4-2\,\zeta(2))-\sum_{n\geq 1}\frac{1}{n^2(2n+1)^2(4n+1)}\\
&=& 48-8\pi-32\log(2)
\end{eqnarray*}
\begin{eqnarray*}&&\int_{0}^{1}K(x)\log^2(1-x)\,dx \stackrel{\text{Taylor}}{=} \frac{\pi}{2}\sum_{n\geq 0}\frac{H_{n+1}^2+H_{n+1}^{(2)}}{n+1}\left[\frac{1}{4^n}\binom{2n}{n}\right]^2 \\
&\stackrel{\text{FL}}{=}& 4+\sum_{n\geq 1}\frac{2}{(2n+1)^2}\left[2\left(\frac{1}{n}+\frac{1}{n+1}\right)^2+4\left(\frac{1}{n}+\frac{1}{n+1}\right)H_{n-1}\right]\\
&=& -8+\frac{4\pi^2}{3}+8\sum_{n\geq 1}\frac{H_{n-1}}{n(n+1)(2n+1)}\\
&=& 48-\frac{4\pi^2}{3}-64\log(2)+32\log^2(2)
\end{eqnarray*}
\begin{eqnarray*}
&&\frac{\pi}{2}\sum_{n\geq 0}\frac{H_{n+1}}{(n+1)^2}\left[\frac{1}{4^n}\binom{2n}{n}\right]^2\\ &\stackrel{\text{Taylor}}{=}& -\int_{0}^{1}\left[4\frac{d}{dx}E(x)+2\,K(x)\right]\log(1-x)\,dx\\
&\stackrel{\text{FL}}{=}&(8-2\pi)+\sum_{n\geq 1}\left(\frac{1}{n}+\frac{1}{n+1}\right)\left[\frac{4}{(2n+1)^2}-\frac{4}{2n+1}+8\int_{0}^{1}\frac{u^{2n+2}}{1+u^2}\,du\right]\\
&=& 16-2\pi-16\log(2)-8\int_{0}^{1}\frac{u^2+(1+u^2)\log(1-u^2)}{1+u^2}\,du\\
&=& 24-32\log(2).
\end{eqnarray*}
By combining these integrals/series we immediately get:
$$ \sum_{n\geq 0}\frac{H_{n+1}^2}{n+1}\left[\frac{1}{4^n}\binom{2n}{n}\right]^2 = 16+\frac{48}{\pi}-\frac{10 \pi }{3}-\frac{128\log(2)}{\pi}+\frac{64\log^2(2)}{\pi}$$
$$ \sum_{n\geq 0}\frac{H_{n+1}^{(2)}}{n+1}\left[\frac{1}{4^n}\binom{2n}{n}\right]^2 = -16+\frac{48}{\pi }+\frac{2 \pi }{3}.$$
Many other identities involving $H_n$ or $\mathscr{H}_n$ can be derived through the FL-method. For instance, we may consider that 
by differentiating Euler's Beta function we have
$$\int_{0}^{1}\frac{\log^2(x)}{\sqrt{x(1-x)}}\,dx = \frac{\pi^3}{3}+4\pi\log^2(2),$$
where the LHS is at the same time
$$ 2\pi + \pi\sum_{n\geq 1}\left[2\left(\frac{1}{2n}+\frac{1}{2n+1}\right)^2+4\left(\frac{1}{2n}+\frac{1}{2n+1}\right)H_{2n-1}\right]\cdot\left[\frac{1}{4^n}\binom{2n}{n}\right]^2 $$
and 
$$ 32\sum_{n\geq 0}(-1)^n\left[\frac{2\mathscr{H}_n^2}{2n+1}-\frac{2\mathscr{H}_n}{(2n+1)^2}+\frac{1}{(2n+1)^3}\right] $$
due to the FL-expansions of $\log^2(x),\frac{1}{\sqrt{x(1-x)}},\frac{\log^2(x)}{\sqrt{x}}$ and $\frac{1}{\sqrt{1-x}}$.\\
In particular the series 
$$\sum_{n\geq 1}\left(\frac{1}{2n}+\frac{1}{2n+1}\right)H_{2n}\left[\frac{1}{4^n}\binom{2n}{n}\right]^2$$ is given by 
an algebraic combination of $1,\log(2),\pi,K$ and $\text{Im}\,\text{Li}_3\left(\tfrac{1+i}{2}\right)$.\\ 
Similarly,
$$ \int_{0}^{1}\frac{-\log^3(x)}{\sqrt{x(1-x)}}\,dx = 2\pi^3\log(2)+8\pi\log^3(2)+12\pi\zeta(3) $$
and the FL-expansion of $\frac{-\log^3(x)}{\sqrt{x}}$ can be computed through the same technique we used for the 
FL-expansion of $\frac{\log^2(x)}{\sqrt{x}}$, namely Rodrigues' formula and recurrence relations.
In $L^2(0,1)$ we have equality between $-\frac{\log^3(x)}{\sqrt{x}}$ and
$$32\sum_{n\geq 0}(-1)^n P_n(2x-1)\left[4\mathscr{H}_n^3 + 2\mathscr{H}_n^{(3)}-\frac{6\mathscr{H}_n^2}{(2n+1)}+\frac{6\mathscr{H}_n}{(2n+1)^2}-\frac{3}{(2n+1)^3}\right] $$
and the FL-expansion of $\log^3(x)$, it is not difficult to find a closed form for $\sum_{n\geq 0}\frac{(-1)^n \mathscr{H}_n^{(3)}}{2n+1}$ through twisted 
hypergeometric series.
Now we have an important remark. According to the \emph{niceness} criterion outlined in \cite{Interplay}, we may denote through $\mathscr{N}$ the class of (twisted) hypergeometric $\phantom{}_{p+1}F_p(x)$ functions such that 
the coefficients of their FL-expansions are the coefficients of the Maclaurin series of a (twisted) $\phantom{}_{q+1}F_q(x)$ function, with $\max(p,q)$ being reasonably small. As already remarked, 
assuming $a+3,b+3\in\frac{1}{4}\N$ and $u\in(-1,1)$ all the following functions belong to $\mathscr{N}$:
$$ \frac{1}{\sqrt{(1+u)^2-4ux}},\quad x^{a}(1-x)^b, \quad \log(x), \quad \log^2(x), \quad \log(1-\sqrt{x}),\quad \frac{\log(x)}{\sqrt{x}},$$
$$ \frac{\log^2(x)}{\sqrt{x}},\quad K(x),\quad E(x),\quad \mathcal{J}(x),\quad \frac{d}{dx}E(x),\quad \frac{\arcsin\sqrt{x}}{\sqrt{x}} $$
and $\mathscr{N}$ is closed with respect to $f(x)\mapsto f(1-x), x\cdot f(1-x), \int f(x)\,dx$ by Bonnet's recursion formulas.
\begin{theorem}\label{li2}$\text{Li}_2(x)\in\mathscr{N}.$\end{theorem}
\begin{proof} It is straightforward to check that $\int_{0}^{1}\text{Li}_2(x)\,dx = \zeta(2)-1$ and that
$$ \int_{0}^{1}\text{Li}_2(x)P_n(2x-1)\,dx = \frac{1}{n^2(n+1)^2} $$
holds for any $n\in\{1,2,3,4,5\}$. This leads to the following conjectural form for the FL-expansion of the dilogarithm:
\begin{equation}\label{conjli2} \text{Li}_2(x) = \left(\zeta(2)-1\right) + \sum_{n\geq 1}\left(\frac{1}{n^2}-\frac{1}{(n+1)^2}\right) P_n(2x-1) \end{equation}
which can be proved to be correct through the following approach: by applying the operator $x\cdot\frac{d}{dx}$ to both sides of (\ref{conjli2}) and exploiting Bonnet's recursion
formulas, one finds $-\log(1-x)$ in the LHS and $1+\sum_{n\geq 1}\left(\frac{1}{n}+\frac{1}{n+1}\right) P_n(2x-1)$ in the RHS. By Rodrigues formula we know 
that the last series is indeed the FL-expansion of $-\log(1-x)$.\qed\end{proof}
\begin{corollary}
Despite their appearance, all the following integrals are very simple to compute:
$$ \int_{0}^{1}\text{Li}_2^2(x)\,dx, \;\; \int_{0}^{1}\text{Li}_2^2(\sqrt{x})\,dx, \;\; \int_{0}^{1}\text{Li}_2(x)\text{Li}_2(1-x)\,dx, \;\; \int_{0}^{1}\text{Li}_2(x) K(x)\,dx, $$
$$ \int_{0}^{1}\frac{\text{Li}_2(x)}{\sqrt{x(1-x)}}\,dx, \quad \int_{0}^{1}\frac{\text{Li}_2(x)}{\sqrt[4]{x(1-x)}}\,dx,\quad \int_{0}^{1}\text{Li}_2(x)\log(x)\log(1-x)\,dx.$$
\end{corollary}
\begin{proof} The integral $\int_{0}^{1}\text{Li}_2^2(x)\,dx$ equals
$$\sum_{m,n\geq 1}\frac{1}{m^2 n^2 (m+n+1)}=\sum_{s\geq 2}\frac{1}{s+1}\sum_{a=1}^{s-1}\frac{1}{a^2(s-a)^2} = \sum_{s\geq 2}\frac{2s H_{s-1}^{(2)}+4 H_{s-1}}{s^3(s+1)} $$
and the RHS can be computed through standard Euler sums. As an alternative, the LHS can be computed by integration by parts, since 
$\int_{0}^{1}\text{Li}_2(x)\log(1-x)\,dx = 3-\zeta(2)-2\zeta(3)$ and $\int \text{Li}_2(x)\,dx = C+x\,\text{Li}_2(x)-(1-x)\log(1-x)-x$. FL-expansions 
provide a more efficient derivation, since by Theorem (\ref{li2}) the LHS equals
$$ \left(\zeta(2)-1\right)^2 + \sum_{n\geq 1}\frac{2n+1}{n^4(n+1)^4} \stackrel{\text{PFD}}{=} 6-2\zeta(2)+\zeta(2)^2-4\zeta(3).$$
Similarly the integral $\int_{0}^{1}\text{Li}_2^2(\sqrt{x})\,dx$ equals 
$$ \sum_{m,n\geq 1}\frac{2}{m^2 n^2 (m+n+2)}=\sum_{s\geq 2}\frac{2}{s+2}\sum_{a=1}^{s-1}\frac{1}{a^2(s-a)^2} = \sum_{s\geq 2}\frac{4s H_{s-1}^{(2)}+8 H_{s-1}}{s^3(s+2)} $$
and by $\text{IBP}$ or Euler sums both sides equal $\frac{25}{8}-\frac{\pi^2}{4}+\frac{\pi ^4}{36}-2\,\zeta(3)$. On the other hand the computation of the FL-expansion of $x\,\text{Li}_2(x)$ provides an equivalent series 
in which harmonic numbers \emph{do not} appear.\\
Through Euler's Beta function, the integral $\int_{0}^{1}\text{Li}_2(x)\text{Li}_2(1-x)\,dx$ can be written as:
$$ \sum_{m,n\geq 1}\frac{1}{m^2 n^2 (m+n+1) \binom{m+n}{n}}=\sum_{s\geq 2}\frac{1}{s^3(s+1)}\sum_{a=1}^{s-1}\frac{s+2a}{a^2 \binom{s}{a}} $$
where the equality $\sum_{a=0}^{n}\binom{n}{a}^{-1}=\frac{n+1}{2^{n+1}}\sum_{k=1}^{n+1}\frac{2^k}{k}$ involving ``generalized'' harmonic numbers follows from 
Euler's Beta function too. The FL-expansion of $\text{Li}_2$ simply grants
$$ \int_{0}^{1}\text{Li}_2(x)\text{Li}_2(1-x)\,dx = \left(\zeta(2)-1\right)^2 + \sum_{n\geq 1}\frac{(-1)^n (2n+1)}{n^4(n+1)^4}\stackrel{\text{PFD}}{=} 6-\frac{2 \pi ^2}{3}+\frac{\pi ^4}{120}.$$
The fourth integral can be written as 
$$ \int_{0}^{1} K(x)\text{Li}_2(x)\,dx = \frac{\pi}{2}\sum_{n\geq 0}\left(\frac{\zeta(2)}{n+1}-\frac{H_{n+1}}{(n+1)^2}\right)\left[\frac{1}{4^n}\binom{2n}{n}\right]^2. $$
Given the FL-expansions of $K(x)$ and $\text{Li}_2(x)$, both sides simply equal
$$ 2\left(\zeta(2)-1\right)+\sum_{n\geq 1}\frac{2}{(2n+1)^2}\left(\frac{1}{n^2}-\frac{1}{(n+1)^2}\right)\stackrel{\text{PFD}}{=} -24+32\log(2)+\frac{\pi^2}{3}. $$
The fifth integral can be represented as
$$ \int_{0}^{1}\frac{\text{Li}_2(x)}{\sqrt{x(1-x)}}\,dx = \sum_{n\geq 0}\left(\frac{2\zeta(2)}{2n+1}-\frac{4H_{n+1/2}}{(2n+1)^2}\right)\left[\frac{1}{4^n}\binom{2n}{n}\right] $$
and $\int_{0}^{\pi/2}\text{Li}_2(\sin^2\theta)\,d\theta=\int_{0}^{\pi/2}\text{Li}_2(\cos^2\theta)\,d\theta$ can be proved to be equal to $\frac{\pi^3}{12}-\pi\log^2(2)$ 
through the dilogarithm reflection formulas. On the other hand, by the FL-expansions of $\text{Li}_2(x)$ and $\frac{1}{\sqrt{x(1-x)}}$ both sides equal
$$ \pi(\zeta(2)-1)+\pi\sum_{n\geq 1}\left(\frac{1}{(2n)^2}-\frac{1}{(2n+1)^2}\right)\left[\frac{1}{4^n}\binom{2n}{n}\right]^2 $$
and the hypergeometric $\phantom{}_{p+1} F_{p}(1)$ values mentioned by the last line have already been computed through the FL machinery. In particular we get
$$ \sum_{n\geq 0}\frac{H_{n+1/2}}{(2n+1)^2}\left[\frac{1}{4^n}\binom{2n}{n}\right]=\frac{\pi}{2}\log^2(2),$$ 
$$\sum_{n\geq 0}\frac{\mathscr{H}_n}{(2n+1)^2}\left[\frac{1}{4^n}\binom{2n}{n}\right]=\frac{3\pi}{4}\log^2(2).$$
The sixth integral equals
$$ \int_{0}^{1}\frac{\text{Li}_2(x)}{\sqrt[4]{x(1-x)}}\,dx = \frac{\Gamma\left(\frac{3}{4}\right)}{\pi \sqrt{2}}\sum_{n\geq 0}\frac{\Gamma\left(n+\frac{1}{4}\right)}{n!}\left(\frac{\zeta(2)}{n+3/4}-\frac{H_{n+3/4}}{(n+3/4)^2}\right)$$ $$= \frac{\pi^{3/2}}{3}\Gamma\left(\tfrac{3}{4}\right)^2 -\frac{\Gamma\left(\frac{3}{4}\right)}{\pi \sqrt{2}}\sum_{n\geq 0}\frac{\Gamma\left(n+\frac{1}{4}\right)H_{n+3/4}}{(n+3/4)^2 n!}$$
and by the FL-expansions of $\text{Li}_2(x)$ and $\left[x(1-x)\right]^{-1/4}$ the LHS can be written as
$$ \frac{4\pi^{3/2}}{\Gamma\left(\frac{1}{4}\right)^2}\left(\zeta(2)-1\right)+\frac{4\pi^{3/2}}{\Gamma\left(\frac{1}{4}\right)^2}\sum_{n\geq 1}\frac{1}{(2n)^2(2n+1)^2}\left[\frac{1}{4^n}\binom{2n}{n}\right]$$
only involving elementary hypergeometric functions. By $\text{PFD}$ the final outcome is
$$ \int_{0}^{1}\frac{\text{Li}_2(x)}{\sqrt[4]{x(1-x)}}\,dx = \frac{2\pi^{3/2}}{\Gamma\left(\frac{1}{4}\right)^2}\left(\frac{5 \pi ^2}{12}+(\pi-4)\left(2+\log(2)\right)-\log^2(2)\right)$$
leading to a closed form for $\phantom{}_4 F_3\left(1,1,1,\tfrac{7}{4};2,2,\tfrac{5}{2};1\right)$. The seventh integral is related to an Euler sum with weight four:
$$ \int_{0}^{1}\text{Li}_2(x)\log(x)\log(1-x)\,dx = \sum_{n\geq 1}\left(\frac{H_{n+1}}{n^2(n+1)^2}-\frac{\zeta(2)-H_{n+1}^{(2)}}{n^2(n+1)}\right). $$
Of course the similarity between the FL expansions of $\text{Li}_2(x)$ and $\log(x)\log(1-x)$ (the latter only involving Legendre polynomials with even degree, 
by symmetry) is not accidental at all, but induced by the dilogarithm reflection formula. In particular we may remove $H_{n+1}$ and $H_{n+1}^{(2)}$ from the RHS 
of the previous line either by noticing that
\begin{eqnarray*}&& 2\int_{0}^{1}\text{Li}_2(x)\log(x)\log(1-x)\,dx\\ &=& \int_{0}^{1}\left(\zeta(2)-\log(x)\log(1-x)\right)\log(x)\log(1-x)\,dx\\ 
&=& \zeta(2)\left.\frac{\partial^2}{\partial a\,\partial b}\left(\frac{\Gamma(a+1)\Gamma(b+1)}{\Gamma(a+b+2)}\right)-\frac{\partial^4}{\partial^2 a\,\partial^2 b}\left(\frac{\Gamma(a+1)\Gamma(b+1)}{\Gamma(a+b+1)}\right)\right|_{(a,b)=(0,0)}\end{eqnarray*}
or by noticing that
\begin{eqnarray*}&&\int_{0}^{1}\text{Li}_2(x)\log(x)\log(1-x)\,dx\\ &\stackrel{\text{FL}}{=}& (\zeta(2)-1)(2-\zeta(2))-\sum_{n\geq 1}\frac{1}{2n^2(2n+1)^2}\left(\frac{1}{(2n)^2}-\frac{1}{(2n+1)^2}\right)\\
&\stackrel{\text{PFD}}{=}& -12+\frac{5\pi^2}{6}+4\,\zeta(3)-\frac{\pi^4}{120}.
\end{eqnarray*}
\qed\end{proof}
\begin{corollary} Since $\int_{0}^{1}\text{Li}_2(x) x^n\,dx = \frac{\zeta(2)}{n+1}-\frac{H_{n+1}}{(n+1)^2}$, the FL-expansion of the dilogarithm function can be efficiently
used in the computation of many twisted hypergeometric series involving $H_n$. 
\end{corollary}
\begin{proof} Since $\int_{0}^{1}K(x)\,dx = \frac{\pi}{2}\sum_{n\geq 0}\frac{1}{(n+1)}\left[\frac{1}{4^n}\binom{2n}{n}\right]^2 = 2$, 
the computation of\\ $\int_{0}^{1}K(x)\text{Li}_2(x)\,dx$ performed in the previous corollary immediately leads to
$$ \sum_{n\geq 0}\frac{H_{n+1}}{(n+1)^2}\left[\frac{1}{4^n}\binom{2n}{n}\right]^2 = \frac{48-64\log 2}{\pi}.$$
Analogously, the FL machinery provides an alternative proof for some instances of Euler's classical result about $\sum_{n\geq 1}\frac{H_n}{n^s}$, 
which is Theorem 2.2 in \cite{Flajolet}. For instance
\begin{eqnarray*} &&\int_{0}^{1}\text{Li}(x)^2\,dx = \sum_{n\geq 0}\int_{0}^{1}\frac{x^{n+1}}{(n+1)^2}\text{Li}_2(x)\,dx\\ &=& \sum_{n\geq 0}\frac{1}{(n+1)^2}\left(\frac{\zeta(2)}{n+1}-\frac{H_{n+1}}{(n+1)^2}\right)=\zeta(2)\zeta(3)-\sum_{n\geq 1}\frac{H_n}{n^4}.\end{eqnarray*}
\qed\end{proof}
\begin{theorem}\label{li3}$\text{Li}_3(x)\in\mathscr{N}.$\end{theorem}
\begin{proof} We start by computing the moments of $\text{Li}_3(x)$ over $(0,1)$:
\begin{eqnarray}\label{momli3}\notag\int_{0}^{1}x^n\text{Li}_3(x)\,dx&=&\sum_{m\geq 1}\frac{1}{m^3(m+n+1)}\\&\stackrel{\text{PFD}}{=}&\frac{\zeta(3)}{n+1}-\frac{\zeta(2)}{(n+1)^2}+\frac{H_{n+1}}{(n+1)^3}.\end{eqnarray}
Since
$$ \sum_{k=0}^{n}\binom{n}{k}\binom{n+k}{k}\frac{(-1)^k}{k+1}=0,\quad \sum_{k=0}^{n}\binom{n}{k}\binom{n+k}{k}\frac{(-1)^k}{(k+1)^2}=\frac{1}{n(n+1)} $$
for any $n\geq 1$ we have
$$ \int_{0}^{1} \text{Li}_3(x)P_n(2x-1)\,dx = (-1)^n\left[u_n - \frac{\zeta(2)}{n(n+1)}\right] $$
and 
\begin{eqnarray*} \text{Li}_3(x) &=& \left(1-\zeta(2)+\zeta(3)\right) + \sum_{n\geq 1}(-1)^n(2n+1)\left[u_n - \frac{\zeta(2)}{n(n+1)}\right]P_n(2x-1)\\
&=& \left(1-\zeta(2)+\zeta(3)\right)+\frac{2\pi^2-15}{8}P_1(2x-1) + \sum_{n\geq 2}a_n P_n(2x-1). \end{eqnarray*}
On the other hand $x\cdot \text{Li}_3(x)-\int \text{Li}_3(x)\,dx = x\cdot\text{Li}_2(x)-x-(1-x)\log(1-x)$ and the FL-expansion 
of the RHS is known. In particular we may find the unknown coefficients $u_n$, associated to the sums
$$ \sum_{k=0}^{n}\binom{n}{k}\binom{n+k}{k}\frac{(-1)^k H_{k+1}}{(k+1)^3} $$
by solving a recurrence relation induced by Bonnet's formulas.\\ Indeed, if $f(x)=\sum a_n P_n(2x-1)$ and $x\cdot f(x)-\int f(x)\,dx=\sum b_n P_n(2x-1)$,\\ we have
$$ b_n = \frac{n-1}{2(2n-1)}a_{n-1}+\frac{1}{2}a_n+\frac{n+2}{2(2n-3)}a_{n+1} $$
for any $n\geq 2$. By computing the FL-expansion of $x\cdot\text{Li}_2(x)-(1-x)\log(1-x)-x$ we get
$$ \frac{n-1}{2(2n-1)}a_{n-1}+\frac{1}{2}a_n+\frac{n+2}{2(2n+3)}a_{n+1} = \frac{2(2n+1)(n^2+n+1)}{(n-1)^2 n^2 (n+1)^2 (n+2)} $$
whose solution $a_n$ can be computed through the WZ method \cite{Wilf}:
\begin{equation}\label{flli3}\frac{(2n+1)\left(n^4+n^3+3n^2+8n+4\right)}{n^3 (n+1)^3 (n+2)^2} + \frac{(2n+1)(-1)^n}{2n(n+1)}\left[4\sum_{k=1}^{n+2}\frac{(-1)^{k+1}}{k^2}-\frac{\pi^2}{3}\right].\end{equation}
The last line provides an explicit expression for the FL-expansion of $\text{Li}_3$ in terms of an alternating quadratic sum (a tail for the classical 
series defining $\eta(2)$) and proves the following identity for any $n\geq 2$:
\begin{eqnarray}\notag&&\sum_{k=0}^{n}\binom{n}{k}\binom{n+k}{k}\frac{(-1)^k H_{k+1}}{(k+1)^3}\\&=&\frac{(-1)^n (n^4+n^3+3n^2+8n+4)}{n^3 (n+1)^3 (n+2)^2} + \frac{2}{n(n+1)}\sum_{k=1}^{n+2}\frac{(-1)^{k+1}}{k^2}.\end{eqnarray}
\qed\end{proof}
\begin{corollary} \begin{eqnarray*} \int_{0}^{1}\frac{\text{Li}_3(x)}{\sqrt{1-x}}\,dx &=& \frac{4}{3}\phantom{}_4 F_3\left(1,1,1,1;2,2,\tfrac{5}{2};1\right)\\
             &=&\sum_{n\geq 0}\left(\frac{\zeta(3)}{n+1}-\frac{\zeta(2)}{(n+1)^2}+\frac{H_{n+1}}{(n+1)^3}\right)\left[\frac{1}{4^n}\binom{2n}{n}\right]\\
             &=&2\zeta(3)-4\zeta(2)+4\zeta(2)\log(2)+\sum_{n\geq 0}\frac{H_{n+1}}{(n+1)^3}\left[\frac{1}{4^n}\binom{2n}{n}\right]\\[-0.2cm]
             &=& \sum_{n\geq 1}\frac{2\cdot 4^n}{(2n+1)n^3 \binom{2n}{n}}\\
             &=& 16 - 2\pi^2 +2\pi^2 \log(2)-7\zeta(3).
            \end{eqnarray*}
\end{corollary}
\begin{proof}
By the FL-expansion of $\text{Li}_3$ the given integral equals
$$ -\frac{13}{18}-\frac{\pi ^2}{2}+6\zeta(3)+2\int_{0}^{1}\frac{x\log x}{x+1}\left(2 x-x^2+2(1-x)\log(1-x)\right)\,dx $$
since 
$$ (-1)^n\left[-\frac{\pi^2}{3}+4\sum_{k=1}^{n+2}\frac{(-1)^{k+1}}{k^2}\right]=\int_{0}^{1} x^n\cdot \frac{4x^2\log(x)}{x+1}\,dx.$$
In particular the FL method allows a conversion into an elementary integral, only involving $\log(x),\log(1-x),x$ and $1\pm x$. It also proves that 
$$\sum_{n\geq 0}\frac{H_{n+1}}{(n+1)^3}\left[\frac{1}{4^n}\binom{2n}{n}\right]=16-\frac{4\pi^2}{3}(1-\log 2)-9\zeta(3).$$
\qed\end{proof}
\begin{corollary} \begin{eqnarray*} \int_{0}^{1}\text{Li}_3(x) K(x)\,dx &=& \frac{\pi}{2}\sum_{n\geq 0}\left(\frac{\zeta(3)}{n+1}-\frac{\zeta(2)}{(n+1)^2}+\frac{H_{n+1}}{(n+1)^3}\right)\left[\frac{1}{4^n}\binom{2n}{n}\right]^2\\
             &=& 2\zeta(3)+\frac{\pi^3}{3}-\frac{4\pi^2}{3}+\frac{\pi}{2}\sum_{n\geq 0}\frac{H_{n+1}}{(n+1)^3}\left[\frac{1}{4^n}\binom{2n}{n}\right]^2\\
             &=& 64-\frac{4 \pi ^2}{3}-96\log(2)+2\zeta(3)-3\pi^3+64K\log(2)\\&-&4\pi\log^2(2)+128\,\text{Im}\,\text{Li}_3\left(\tfrac{1+i}{2}\right).
            \end{eqnarray*}
\end{corollary}
\begin{proof} By the FL-expansion of $\text{Li}_3$ the given integral equals
\begin{eqnarray*}&&\frac{451}{9}-\frac{22}{9}\pi^2-\frac{352}{9}\log(2)+2\zeta(3)\\&+&\int_{0}^{1}\frac{4x\log(x)}{x+1}\left(3 x-4\sqrt{x}\,\text{arctanh}\sqrt{x}-(x+1)\log(1-x)\right)\,dx \end{eqnarray*}
which after some elementary manipulations becomes 
$$64-\frac{4 \pi ^2}{3}-96\log(2)+2\zeta(3)+32 \int_{0}^{1}\log\left(\frac{1-x}{1+x}\right)\frac{\log(x)}{1+x^2}\,dx.$$
The functional identities for $\text{Li}_3$ allow to turn the last line into 
$$ 64-\frac{4 \pi ^2}{3}-96\log(2)+2\zeta(3)-3\pi^3+64K\log(2)-4\pi\log^2(2)+128\,\text{Im}\,\text{Li}_3\left(\tfrac{1+i}{2}\right) $$ 
proving the claim and the equality
\begin{eqnarray*} &&\frac{\pi}{2}\sum_{n\geq 0}\frac{H_{n+1}}{(n+1)^3}\left[\frac{1}{4^n}\binom{2n}{n}\right]^2\\ &=& 64 - \frac{10\pi^3}{3}+ 128\,\text{Im}\,\text{Li}_3\left(\tfrac{1+i}{2}\right) - 96\log(2) + 
 64 K \log(2) - 4\pi\log^2(2). \end{eqnarray*}
\qed\end{proof}

\section{Further remarks}
In the previous section we have seen that FL-expansions reduce the problem of evaluating 
$$ \phantom{}_{\mu+3} F_{\mu+2}\left(\tfrac{1}{2},\tfrac{1}{2},1,1,\ldots;2,2,2,\ldots;1\right)=\sum_{n\geq 0}\frac{1}{(n+1)^{\mu+2}}\left[\frac{1}{4^n}\binom{2n}{n}\right]^2,\quad \mu\in\mathbb{N} $$
to the problem of evaluating 
\begin{eqnarray*} \int_{0}^{\pi/4}\log^\mu(\cos\theta)\,d\theta &=& \sum_{n\geq 0}\frac{1}{2^n(2n+1)^{\mu+1}}\left[\frac{1}{4^n}\binom{2n}{n}\right]\\ &=& \phantom{}_{\mu+2} F_{\mu+1}\left(\tfrac{1}{2},\tfrac{1}{2},\ldots;\tfrac{3}{2},\tfrac{3}{2},\ldots;\tfrac{1}{2}\right).\end{eqnarray*}
We have to stress that the LHS of the last line is related to twisted hypergeometric series by a classical approach, which is to consider some combinatorial/hypergeometric 
identities like the Chu-Vandermonde or the Dougall-Dixon ones, then differentiate them with respect to one/some of their parameters. This approach allowed 
Wang and Jia in \cite{Wang} to prove a large number of identies involving generalized harmonic numbers and binomial transforms. In our case we may consider that
\begin{eqnarray}\label{hyp1} \notag&&\int_{0}^{\pi/4}\left(\cos\theta\right)^a\,d\theta = \frac{\sqrt{\pi}\,\Gamma\left(\frac{a}{2}+\frac{1}{2}\right)}{2\,\Gamma\left(\frac{a}{2}+1\right)}-\frac{1}{(a+1) 2^{\frac{a+1}{2}}}\,\phantom{}_2 F_1\left(\tfrac{1}{2},\tfrac{a+1}{2};\tfrac{a+3}{2};\tfrac{1}{2}\right)\\
 &=& \frac{\sqrt{\pi}\,\Gamma\left(\frac{a}{2}+\frac{1}{2}\right)}{2\,\Gamma\left(\frac{a}{2}+1\right)}-\frac{1}{(a+1) 2^{\frac{a+1}{2}}}\sum_{n\geq 0}\frac{1}{8^n}\binom{2n}{n}\left[1-\frac{2n}{2n+a+1}\right]
\end{eqnarray}
follows from the cosine duplication formula, whereas
\begin{equation}\label{hyp2} \int_{0}^{\pi/4}\left(\cos\theta\right)^a\,d\theta = \phantom{}_2 F_1\left(\tfrac{1}{2},\tfrac{a+2}{2};\tfrac{3}{2};-1\right)=\sum_{n\geq 0}\frac{(-1)^n }{(2n+1)n!}\cdot\frac{\Gamma\left(n+\frac{a}{2}+1\right)}{\Gamma\left(\frac{a}{2}+1\right)}\end{equation}
follows from the substitution $\theta\mapsto\arctan u$. By differentiating the right-hand sides of (\ref{hyp1}) and (\ref{hyp2}) multiple times with respect to $a$, 
then performing an evaluation at $a=0$, we get that the values of $\phantom{}_{\mu+2} F_{\mu+1}\left(\tfrac{1}{2},\tfrac{1}{2},\ldots;\tfrac{3}{2},\tfrac{3}{2},\ldots;\tfrac{1}{2}\right)$ or 
$\int_{0}^{\pi/4}\log^\mu(\cos\theta)\,d\theta$ provide explicit forms both for $\sum_{n\geq 0}\frac{1}{(n+1)^{\mu+2}}\left[\frac{1}{4^n}\binom{2n}{n}\right]^2$ and for the twisted sums
\begin{eqnarray}\label{hyp3}\notag\sum_{n\geq 0}\frac{(-1)^n}{(2n+1)}H_n,\qquad \sum_{n\geq 0}\frac{(-1)^n}{(2n+1)}\left(H_n^2-H_n^{(2)}\right),\\ \sum_{n\geq 0}\frac{(-1)^n}{(2n+1)}\left(H_n^3-3 H_n H_n^{(2)}+2H_n^{(3)}\right), \end{eqnarray}
the latter being clearly related to $\sum_{n\geq 0}\frac{(-1)^n \mathscr{H}_n^3}{2n+1}$, which has been crucial for the evaluation of
$\sum_{n\geq 0}\frac{1}{(2n+1)^3}\left[\frac{1}{4^n}\binom{2n}{n}\right]^2 =\phantom{}_5 F_4\left(\tfrac{1}{2},\tfrac{1}{2},\tfrac{1}{2},\tfrac{1}{2},\tfrac{1}{2};1,\tfrac{3}{2},\tfrac{3}{2},\tfrac{3}{2};1\right).$
By exploiting the identity $\log\left(2\sin\frac{x}{2}\right)=\log(1-e^{ix})+\frac{i}{2}(\pi-x)$, Shalev has shown in \cite{Web} that
\begin{eqnarray*} \int_{0}^{\pi/4} \log^3(\sin\theta)\,d\theta &=& 3\,\text{Im}\,\text{Li}_4(1-i)-\frac{25\pi^3}{256}\log(2)+\frac{3}{2}\log(2)\,\text{Im}\,\text{Li}_3\left(\tfrac{1+i}{2}\right)\\
&& -\frac{3K}{8}\log^2(2)-\frac{17\pi}{64}\log^3(2)-\frac{3\pi}{8}\zeta(3)+\frac{3}{4}\beta(4)
\end{eqnarray*}
and by differentiating the Euler Beta function three times it follows that
\begin{eqnarray}\label{intlog3}\notag \int_{0}^{\pi/4} \log^3(\cos\theta)\,d\theta &=& -\frac{7\pi^3}{256}\log(2)-\frac{15\pi}{64}\log^3(2)+\frac{3K}{8}\log^2(2)-\frac{3\pi}{8}\zeta(3)\\
&& -\frac{3}{4}\beta(4)-\frac{3}{2}\log(2)\,\text{Im}\,\text{Li}_3\left(\tfrac{1+i}{2}\right)-3\,\text{Im}\,\text{Li}_4(1-i)
\end{eqnarray}
which provides closed forms for
$$ \sum_{n\geq 0}\frac{1}{(n+1)^5}\left[\frac{1}{4^n}\binom{2n}{n}\right]^2,\quad \sum_{n\geq 0}\frac{1}{2^n(2n+1)^4}\left[\frac{1}{4^n}\binom{2n}{n}\right],$$
$$ \sum_{n\geq 0}\frac{(-1)^n}{(2n+1)}\left(H_n^3-3 H_n H_n^{(2)}+2H_n^{(3)}\right).$$
On the other hand the class of hypergeometric functions of the\\ $\phantom{}_{p+1} F_p\left(\tfrac{1}{2},\tfrac{1}{2},\ldots;\tfrac{3}{2},\tfrac{3}{2},\ldots;x\right)$ kind is very 
peculiar, since all these functions have a closed form in terms of polylogarithms. As shown by Shalev and Reshetnikov in \cite{Reshe} we have
\begin{equation}\label{resh} \phantom{}_4 F_3\left(\tfrac{1}{2},\tfrac{1}{2},\tfrac{1}{2},\tfrac{1}{2};\tfrac{3}{2},\tfrac{3}{2},\tfrac{3}{2};\sin^2\theta\right) = \frac{\theta^3}{6\sin\theta}+\frac{1}{2\sin\theta}\,\text{Im}\,\text{Li}_3\left(2\sin^2\theta+2i\sin\theta\cos\theta\right)\end{equation}
for any $\theta\in\left(0,\frac{\pi}{2}\right)$, hence
\begin{eqnarray*}&&\phantom{}_5 F_4\left(\tfrac{1}{2},\tfrac{1}{2},\tfrac{1}{2},\tfrac{1}{2},\tfrac{1}{2};1,\tfrac{3}{2},\tfrac{3}{2},\tfrac{3}{2};1\right)= \sum_{n\geq 0}\frac{1}{(2n+1)^3}\left[\frac{1}{4^n}\binom{2n}{n}\right]^2 \\
&=& \frac{2}{\pi}\int_{0}^{\pi/2}\frac{\theta^3}{6\sin\theta}+\frac{1}{2\sin\theta}\,\text{Im}\,\text{Li}_3\left(2\sin^2\theta+2i\sin\theta\cos\theta\right)\,d\theta\\
&=& \frac{\pi K}{2}-\frac{4}{\pi}\beta(4)+\frac{1}{\pi}\,\text{Im}\int_{0}^{\pi/2}\text{Li}_3\left(2\sin^2\theta+2i\sin\theta\cos\theta\right)\,\frac{d\theta}{\sin\theta}\\
&=& \frac{\pi K}{2}-\frac{4}{\pi}\beta(4)+\frac{1}{\pi}\,\text{Im}\int_{0}^{+\infty}\text{Li}_3\left(\frac{2}{1-it}\right)\frac{dt}{\sqrt{t^2+1}}\,dt\\
&=& \frac{\pi K}{2}-\frac{4}{\pi}\beta(4)+\frac{1}{\pi}\,\text{Re}\int_{1}^{2}\frac{\text{Li}_3(z)}{z\sqrt{z-1}}\,dz\\
&\stackrel{\text{Li}_3}{=}& \frac{\pi K}{2}-\frac{4}{\pi}\beta(4)+\frac{2}{\pi}\int_{0}^{\pi/4}\text{Li}_3(\cos^2\theta)+\frac{4}{3}\log^3(\cos\theta)-\frac{2\pi^2}{3}\log(\cos\theta)\,d\theta
\end{eqnarray*}
and the evaluation of $\sum_{n\geq 0}\frac{(-1)^n \mathscr{H}_n^3}{(2n+1)}$ or $\sum_{n\geq 0}\frac{(-1)^n \mathscr{H}_n^2}{(2n+1)^2}$ reduces to the computation of
$$ \int_{0}^{1/2}\frac{\text{Li}_3(1-z)}{\sqrt{z(1-z)}}\,dz\quad\text{or}\quad \int_{1/2}^{1}\frac{\text{Li}_2(x)\arcsin\sqrt{x}}{x}\,dx. $$
By using Shalev's form of (\ref{resh}) we have
\begin{equation}\label{shalev} \int_{0}^{1}\frac{\log^2(x)}{\sqrt{x(1- x\sin^2\theta )}}\,dx = \frac{8}{\sin\theta}\left[\frac{\theta^3}{3}-\text{Im}\,\text{Li}_3\left(1-e^{2i\theta}\right)\right]\end{equation}
hence by Fubini's theorem
\begin{eqnarray}\label{complex1}\notag \sum_{n\geq 0}\frac{1}{(2n+1)^3}\left[\frac{1}{4^n}\binom{2n}{n}\right]^2 &=& \frac{\pi K}{2}-\frac{4}{\pi}\beta(4)+\frac{2}{\pi}\,\text{Im}\int_{\gamma}\frac{\text{Li}_3(1-z^2)}{1-z^2}\,dz\\ &=& \frac{\pi K}{2}-\frac{4}{\pi}\beta(4)+\frac{2}{\pi}\,\text{Im}\int_{z_0}^{i}\frac{\text{Li}_3(1-z^2)}{1-z^2}\,dz\end{eqnarray}
where $\gamma$ is the positive-oriented curve $\left\{(\cos\theta,\sin\theta),\theta\in\left(0,\frac{\pi}{2}\right)\right\}$ and $z_0$ is any point $\in(-1,1)$. 
By the generating function for Legendre polynomials and the previously computed $\int_{0}^{1}\frac{\log^2(x)}{\sqrt{x}}P_n(2x-1)\,dx$ we also have
$$\forall\theta\in\left(0,\tfrac{\pi}{2}\right),\qquad \frac{1}{4}\left[\frac{\theta^3}{3}-\text{Im}\,\text{Li}_3(1-e^{2i\theta})\right]=\sum_{n\geq 0}\left(\tan\tfrac{\theta}{2}\right)^{2n+1}(-1)^n J_n,$$
$$\forall t\in(0,1),\qquad \frac{2}{3}\arctan^3(t)-\frac{1}{4}\text{Im}\,\text{Li}_3\left(\frac{4it}{(i+t)^2}\right)=\sum_{n\geq 0}t^{2n+1} (-1)^n J_n,$$
$$ J_n = \frac{2\mathscr{H}_n^2}{2n+1}-\frac{2\mathscr{H}_n}{(2n+1)^2}+\frac{1}{(2n+1)^3}.$$
By the FL-expansion of $\text{Li}_3$, the computation of $\sum_{n\geq 0}\frac{1}{(2n+1)^3}\left[\frac{1}{4^n}\binom{2n}{n}\right]^2$ is equivalent to 
the computation of 
$$\sum_{n\geq 1}\frac{(-1)^n \sum_{k=1}^{n}\frac{(-1)^k}{k^2}}{n(n+1)}\left[\frac{1}{4^n}\binom{2n}{n}\right]^2, $$
and by solving, through the inverse Laplace transform, the moment problem
$$ \int_{0}^{1} x^n f(x)\,dx = (-1)^n \sum_{k=1}^{n}\frac{(-1)^k}{k^2} $$
we get that the computation of the previous series is equivalent to the computation of the integral
$$ \int_{0}^{1}\frac{\log x}{1+x}\sum_{n\geq 0}x^n\left(\frac{1}{n}-\frac{1}{n+1}\right)\left[\frac{1}{4^n}\binom{2n}{n}\right]^2\,dx $$
or the computation of the integrals $\int_{0}^{1}K(x)\log(1+x)\,dx$ and\\ $\int_{0}^{1}K(x)\log(x)\log(1+x)\,dx$.

$$\phantom{}$$
About other classical approaches, especially related to (weighted) sums of squares, an interesting consequence of the Lagrange identity
$$\left(\sum_{k=1}^{n}x_k^2\right)\left(\sum_{k=1}^{n}y_k^2\right)-\left(\sum_{k=1}^{n}x_k y_k\right)^2 = \sum_{1\leq j < k \leq n}(x_j y_k - x_k y_j)^2$$
is the following one:
\begin{lemma}[Knopp]
If $\{a_n\}_{n\geq 0}\in\ell^2(\R)$ is an eventually monotonic sequence, the following series are convergent
$$s\stackrel{\text{def}}{=}\sum_{n=0}^\infty(-1)^na_n,\qquad \delta_k\stackrel{\text{def}}{=}\sum_{n=0}^\infty a_na_{n+k}, \qquad \Delta\stackrel{\text{def}}{=}\sum_{k=1}^\infty(-1)^{k-1}\delta_k $$
and we have:
$$ \sum_{n\geq 0} a_n^2 = s^2+2\Delta. $$
\end{lemma}
Knopp's lemma is tipically used for proving $\zeta(2)=\frac{\pi^2}{6}$ through elementary manipulations, once $a_n$ has been defined as $\frac{1}{n+1}$. 
On the other hand, it has interesting consequences for different hypergeometric sequences, too. Let us define $a_n$ as 
$$ a_n = \frac{1}{1-2n}\left[\frac{1}{4^n}\binom{2n}{n}\right] = [x^n]\sqrt{1-x}.$$
For any $k\in\mathbb{N}^+$ we have
\begin{eqnarray*} \delta_k=\sum_{n\geq 0}a_n a_{n+k}&=&\frac{1}{2\pi}\int_{0}^{2\pi} \sqrt{1-e^{i\theta}}\sqrt{1-e^{-i\theta}} e^{-ki\theta}\,d\theta\\&=&\frac{1}{\pi}\int_{0}^{2\pi}\left|\sin\tfrac{\theta}{2}\right|\cos(k\theta)\,d\theta = \frac{4}{\pi(1-4k^2)}\end{eqnarray*}
and $s=\sqrt{2}$. By partial fraction decomposition we have $\Delta= \frac{2}{\pi}-1$, hence Knopp's lemma immediately leads to an alternative proof of
\begin{equation}\label{quattrosupi} \sum_{n\geq 0}\frac{1}{(2n-1)^2}\left[\frac{1}{4^n}\binom{2n}{n}\right]^2 = \phantom{}_2 F_1 \left(-\tfrac{1}{2},-\tfrac{1}{2};1;1\right)= \frac{4}{\pi}. \end{equation}
In a similar fashion, if we define $a_n$ as $\frac{1}{n+1}\left[\frac{1}{4^n}\binom{2n}{n}\right]=[x^n]\frac{2}{1+\sqrt{1-x}}$, we have 
$$ \delta_k = \sum_{n\geq 0}a_n a_{n+k} = \frac{4}{2k-1}\cdot\frac{1}{4^{k}}\binom{2k}{k}-\frac{16}{\pi(4k^2-1)}$$
and $s=2(\sqrt{2}-1)$, $\Delta=-8+4\sqrt{2}+\frac{8}{\pi}$, from which
$$ \sum_{n\geq 0}\frac{1}{(n+1)^2}\left[\frac{1}{4^n}\binom{2n}{n}\right]^2 = \frac{16}{\pi}-4.$$
The last identity can also be derived from (\ref{quattrosupi}) by reindexing, but of course the interesting part of this approach 
does not rely on the mere computations, but on the fact that it provides an alternative way for simplifying the structure of a hypergeometric series.
On the other hand the class of hypergeometric terms $a_n$ such that $\delta_k$ has a reasonably concise form appears to be pretty narrow:
we plan to perform a thorough analysis of the full potential of Knopp's lemma in future works.\\
$\phantom{}$\\
A third remark is about a symmetry trick. Campbell recently showed (in \cite{Campbell2}) how to compute $\sum_{n\geq 1}\left[\frac{1}{4^n}\binom{2n}{n}\right]^2 \frac{H_n}{2^n}$
(and similar series) from an integral transform, but a slightly different route proves that such series is \emph{clearly} given by a combination 
of $\pi$ and values of the logarithm and the complete elliptic integral of the first kind. Indeed, by directly exploiting the Maclaurin series of $K(x)$ we have:
\begin{eqnarray*}
 &&\sum_{n\geq 1}\left[\frac{1}{4^n}\binom{2n}{n}\right]^2 \frac{H_n}{2^n} = \frac{2}{\pi}\int_{0}^{1}\frac{K(u/2)-K(1/2)}{u-1}\,du 
\\ & \stackrel{\text{Fubini}}{=}& \frac{2}{\pi}\int_{0}^{\pi/2}\int_{0}^{1}\left(\frac{1}{\sqrt{1-\frac{u}{2}\sin^2\theta}}-\frac{1}{\sqrt{1-\frac{1}{2}\sin^2\theta}}\right)\frac{du}{u-1}\,d\theta
\end{eqnarray*}
and the evaluation of the RHS boils down to the evaluation of $\int_{0}^{\pi}\frac{\log(3-\cos\varphi)}{\sqrt{3-\cos\varphi}}\,d\varphi$ or $\int_{0}^{+\infty}\frac{\log(t)}{\sqrt{(t^2+1)(2t^2+1)}}\,dt$.
On the other hand, if $p(x)$ is a real, quadratic and palindromic polynomial, attaining positive values over $\mathbb{R}^+$, then $\int_{0}^{+\infty}\frac{\log(x)}{p(x)}\,dx$ equals zero  
as a consequence of the substitution $x\mapsto \frac{1}{x}$. By the very same reason, if $q(x)$ is a real, biquadratic and palindromic polynomial, attaining positive 
values over $\mathbb{R}^+$, then $\int_{0}^{+\infty}\frac{\log(x)}{\sqrt{q(x)}}\,dx = 0.$ In particular the rapidly convergent twisted hypergeometric series 
$\sum_{n\geq 1}\left[\frac{1}{4^n}\binom{2n}{n}\right]^2 \frac{H_n}{2^n}$ has a closed form in terms of $\pi,\log(2)$ and $K\left(\frac{1}{2}\right)=\frac{1}{4\sqrt{\pi}}\Gamma\left(\frac{1}{4}\right)^2$ only. 
Given the relation between $\text{{\bf K}}(k_1)$ and $\text{{\bf K}}(k_4)$, this principle allows us to state that
$$ \sum_{n\geq 1}\left[\frac{1}{4^n}\binom{2n}{n}\right]^2\frac{H_n}{(\sqrt{2}+1)^{4n}} $$
has a closed form in terms of $\sqrt{2},\pi,\log(2),\log(1+\sqrt{2})$ and the lemniscate constant \emph{without actually performing any computation}.
The same applies to 
$$ \sum_{n\geq 1}\left[\frac{1}{4^n}\binom{2n}{n}\right]^2\frac{H_n}{(n+1)2^n}, $$
which depends on $\pi,\log(2),\Gamma\left(\frac{1}{4}\right)$ and its reciprocal, since $K\left(\frac{1}{2}\right)$ and $E\left(\frac{1}{2}\right)$ 
are related via Legendre's identity. Analogously, a closed form of
$$ \sum_{n\geq 1}\left[\frac{1}{4^n}\binom{2n}{n}\right]^2\frac{H_n}{(n+1)(\sqrt{2}+1)^{4n}} $$
only depends on $\pi,\log(2),\log(\sqrt{2}+1),\Gamma\left(\frac{1}{4}\right)$ and its reciprocal. This leaves us with an interesting topic 
for future investigations, partially covered by the first section and the last paragraph:
\begin{center}\emph{
Is it possible to foresee, through symmetry tricks and/or\\ reflection/transformation formulas for hypergeometric functions\\ and polylogarithms, 
the constants involved in the closed form\\ for a (twisted) hypergeometric series?}
\end{center}
The last remark is about the operator 
 $$ \mathscr{C\!D}: f(x)\mapsto \frac{1}{x}\cdot \mathcal{L}^{-1}\left((2n+1)\int_{0}^{1}f(x)P_n(2x-1)\,dx\right)(-\log x) $$
defined on the space of functions in $L^2(0,1)$, such that the FL-expansion has an inverse Laplace transform. This operator 
is self-adjoint by construction, and the complexity of $\mathscr{C\!D}(f)$ measures the effectiveness of the FL-method in providing 
a closed form to $\int_{0}^{1}f(x)g(x)\,dx$. Of course if $f\in L^2(0,1)$ is such that $f(x)=f(1-x)$, all the coefficients
of the FL-expansion of $f(x)$ with odd index are zero, hence it is practical to introduce the similar operator
$$ \mathscr{C\!D}_0: f(x)\mapsto \frac{1}{x}\cdot \mathcal{L}^{-1}\left((4n+1)\int_{0}^{1}f(x)P_{2n}(2x-1)\,dx\right)(-\log x) $$
only acting on the coefficients of the FL-expansion with an even index. We introduce a compendium of some functions we have dealt with through this article, 
together with their $\mathscr{C\!D}$ or $\mathscr{C\!D}_0$-transforms:
%\begin{center}
%\begin{TAB}(r,0.6cm,0.6cm)[5pt]{|c|c|}{|c|c|c|c|c|c|c|c|}% (rows,min,max)[tabcolsep]{columns}{rows}
%$f(x)$ & $(\mathscr{C\!D} f)(x)$ or $(\mathscr{C\!D}_0 f)(x)$    \\
%$\sqrt{1-x}$ & $\frac{1}{4}\left(\sqrt{x}-\frac{1}{x\sqrt{x}}\right)$  \\
%$K(x)$ & $\frac{1}{\sqrt{x}}$ \\
%$E(x)$ & $\frac{1}{4}\left(-\frac{1}{x\sqrt{x}}+\frac{2}{\sqrt{x}}-\sqrt{x}\right)$ \\ 
%$-\log(1-x)$ & $1+\frac{1}{x}$ \\
%$\log^2(1+\frac{1}{\sqrt{1-x}})$ & $1$ \\
%$\log(x)\log(1-x)$ & $\log\left(x\right)\frac{1-\sqrt{x}}{2x}$ \\
%$\text{Li}_2(x)$ & $\left(1-\frac{1}{x}\right)\log(x)$\\
%\end{TAB}
%\end{center}
\renewcommand\arraystretch{1.6}
\begin{center}
\begin{tabular}{|c|c|}
\hline $f(x)$ & $(\mathscr{C\!D} f)(x)$ or $(\mathscr{C\!D}_0 f)(x)$  \\
\hline $\sqrt{1-x}$ & $\frac{1}{4}\left(\sqrt{x}-\frac{1}{x\sqrt{x}}\right)$  \\
\hline $K(x)$ & $\frac{1}{\sqrt{x}}$ \\
\hline $E(x)$ & $\frac{1}{4}\left(-\frac{1}{x\sqrt{x}}+\frac{2}{\sqrt{x}}-\sqrt{x}\right)$ \\ 
\hline $-\log(1-x)$ & $1+\frac{1}{x}$ \\
\hline $\log^2(1+\frac{1}{\sqrt{1-x}})$ & $1$ \\
\hline $\log(x)\log(1-x)$ & $\log\left(x\right)\frac{1-\sqrt{x}}{2x}$ \\
\hline $\text{Li}_2(x)$ & $\left(1-\frac{1}{x}\right)\log(x)$\\
\hline
\end{tabular}
\end{center}

\section{Conclusions}
Continuing from our previous work, we performed a deeper investigation about the relations between (twisted) hypergeometric sums and Euler sums / polylogarithms. 
The geometry of hypergeometric series with half-integer or quarter-integer parameters, from the contiguity relations to the Gauss-Kummer transformations, turns out to be heavily related to the geometry of Euler sums
and the arrangements of hyperplanes (the \emph{shuffle relations} mentioned by Flajolet and Salvy in \cite{Flajolet} are just birational maps of hypercubes into 
themselves). The accumulated evidence brings us to state that in many practical cases the FL machinery is extremely effective in simplifying the structure 
of a (twisted) hypergeometric series, especially if combined with the computation of residues for suitable kernels involving $\psi(z)=\frac{d}{dz}\log\Gamma(z)$, 
or with differentiation of combinatorial/hypergeometric identities (Gauss, Chu-Vandermonde, Dougall-Dixon). We reduced a large number of $\phantom{}_{p+1} F_{p}(\ldots;\pm 1)$ (and \emph{twisted}) values to standard mathematical 
constants, out-performing the capabilities of many computer algebra systems. In a forthcoming future, it might be the case to devise a fully automated conversion 
procedure for $\phantom{}_{p+1} F_{p}(\ldots;\pm 1)$ values, driven by FL-expansions and Euler sums, to be coupled with the classical 
Wilf-Zeilberger method \cite{Wilf},\cite{Guillera}. That would largely increase the potential of computer algebra systems in the symbolic evaluation of 
hypergeometric series and twisted hypergeometric series.

\section{Compliance with ethical standards}
The authors declare that they have no conflict of interest.
%

%

%\begin{acknowledgements}
%If you'd like to thank anyone, place your comments here
%and remove the percent signs.
%\end{acknowledgements}

% BibTeX users please use one of
%\bibliographystyle{spbasic}      % basic style, author-year citations
%\bibliographystyle{spmpsci}      % mathematics and physical sciences
%\bibliographystyle{spphys}       % APS-like style for physics
%\bibliography{}   % name your BibTeX data base

\begin{thebibliography}{24}
\bibitem[1]{Andrews} G.E. Andrews, R. Askey, R. Roy, \emph{Special Functions}, Cambridge Univ. Press (2000)
\bibitem[2]{Bauer} G. Bauer, \emph{Von den Coefficienten der Reihen von Kugelfunctionen einer Variabeln}, J. Reine Angew. Math. 56, 101–121 (1859) 
\bibitem[3]{Borwein} D. Borwein, J. Borwein, R. Girgensohn, \emph{Explicit evaluation of Euler sums}, Proc. Edinburgh Math. Soc., 38 (1995)
\bibitem[4]{AGM} J. Borwein, P. Borwein, \emph{$\pi$ and the AGM. A Study in Analytic Number Theory and Computational Complexity}, Wiley, New York (1987) 
\bibitem[5]{Campbell1} J. M. Campbell, \emph{Ramanujan-like series for $\frac{1}{\pi}$ involving harmonic numbers}, HAL:01364815 (2017)
\bibitem[6]{Campbell2} J. M. Campbell, \emph{New series involving harmonic numbers and squared central binomial coefficients}, HAL:01774708 (2018)
\bibitem[7]{Parbelos} J. M. Campbell, J. D'Aurizio, J. Sondow, \emph{Hypergeometry of the Parbelos}, to appear in The American Mathematical Monthly (2018)
\bibitem[8]{Interplay} J. M. Campbell, J. D'Aurizio, J. Sondow, \emph{On the interplay between hypergeometric functions, complete elliptic integrals and Fourier-Legendre series expansions}, arXiv:1710.03221 (2017)
\bibitem[9]{Ditrani} J. D'Aurizio, S. Ditrani, \emph{Surprising identities for the hypergeometric ${}_4 F_3$ function}, Bollettino UMI (2017)
\bibitem[10]{DeDoelder} P.J. De Doelder, \emph{On some series containing $\psi(x)-\psi(y)$ and $(\psi(x)-\psi(y))^2$ for certain values of $x$ and $y$}, J. Comput. Appl. Math. 37 (1–3) (1991)
\bibitem[11]{Flajolet} F. Flajolet and B. Salvy, \emph{Euler sums and contour integral representations}, Exp. Math., 7 (1998)
\bibitem[12]{Guillera} J. Guillera, \emph{On WZ-pairs which prove Ramanujan series}, Ramanujan J. 22:249 (2010)
\bibitem[13]{Levrie} P. Levrie, \emph{Using Fourier–Legendre expansions to derive series for $\frac{1}{\pi}$ and $\frac{1}{\pi^2}$}, Ramanujan J. 22:221 (2010)
\bibitem[14]{Viola} R. Marcovecchio, C. Viola, \emph{Irrationality and nonquadraticity measures for logarithms of algebraic numbers}, J. Australian Math. Soc. 92:2 (2012)
\bibitem[15]{Reshe} V. Reshetnikov, N. Shalev \emph{On a closed form for} $\int_{0}^{1}\frac{\log^2(x)\,dx}{\sqrt{x(a-bx)}}$, \url{https://math.stackexchange.com/questions/1438381} (2015)
\bibitem[16]{Web} N. Shalev, \emph{On a closed form for} $\sum_{n\geq 1}\frac{\binom{2n}{n}}{16^n (2n+1)^4}$, \url{https://math.stackexchange.com/questions/915054} (2014)
\bibitem[17]{Sita} R. Sitaramachandra Rao, \emph{A formula of S. Ramanujan}, J. Number Theory 25 (1) (1987)
\bibitem[18]{Vidunas} R. Vidunas, \emph{A generalization of Clausen’s identity}, Ramanujan J. 26:133 (2011) 
\bibitem[19]{XYZ} Ce Xu, Yingyue Yang, Jianwen Zhang, \emph{Explicit evaluation of quadratic Euler sums}, Int. J. Number Theory 13, 655 (2017)
\bibitem[20]{Zagier} D. Zagier, \emph{The Dilogarithm Function}, Frontiers in Number Theory, Physics, and Geometry II. Springer, Berlin, Heidelberg (2007)
\bibitem[21]{Zudilin} W. Zudilin, \emph{Two hypergeometric tales and a new irrationality measure of $\zeta(2)$}, Ann. Math. Québec 38:101 (2014)
\bibitem[22]{Zeng} De-Yin Zheng, \emph{Further summation formulae related to generalized harmonic numbers}, J. Math. Anal. Appl. 335 (2007)
\bibitem[23]{Wang} Weiping Wang, Cangzhi Jia, \emph{Harmonic number identities via the Newton–Andrews method}, Ramanujan J. 35:263–285 (2014)
\bibitem[24]{Wilf} H.S. Wilf, D. Zeilberger, \emph{Rational functions certify combinatorial identities}, J. Am. Math. Soc. 3, 147–158 (1990)
\end{thebibliography}

% Non-BibTeX users please use

\end{document}